\documentclass[psamsfonts]{amsart}

\usepackage{amssymb,amsfonts}
\usepackage{verbatim}
\usepackage[all,arc]{xy}
\usepackage{enumerate}
\usepackage{mathrsfs}
\usepackage[hidelinks]{hyperref}
\usepackage{graphicx}
\usepackage{mathtools}
\usepackage[dvipsnames]{xcolor}
\usepackage{fixme}
\usepackage{physics}
\fxsetup{status=draft} 
\usepackage{svg}
\usepackage{tikz}
\usepackage[toc,page]{appendix}
\usepackage{caption}
\usepackage{subcaption}
\usepackage{multicol}

\usepackage[a4paper, margin=1in]{geometry}

\def \R {\mathbb R }
\def \N {\mathbb N }

\def \Q {\mathbb Q}
\def \Z {\mathbb Z}

\def \C{\mathbb C }
\def \P{\mathbb P }

\renewcommand{\angle}[1]{ \langle #1 \rangle}
\newcommand{\set}[1]{\left\lbrace #1 \right\rbrace}
\renewcommand{\phi}{\varphi}
\newcommand\restr[2]{{\left.\kern-\nulldelimiterspace #1 \vphantom{\big|} \right|_{#2} }}
\newtheorem{thm}{Theorem}[section]

\newtheorem{prop}[thm]{Proposition}
\newtheorem{lem}[thm]{Lemma}

\newtheorem{claim}[thm]{Claim}

\theoremstyle{definition}
\newtheorem{defn}[thm]{Definition}

\newtheorem{notn}[thm]{Notation}

\newtheorem{rem}[thm]{Remark}

\theoremstyle{remark}

\makeatletter
\let\c@equation\c@thm
\makeatother
\numberwithin{equation}{section}

\usepackage{graphicx}
\newcommand{\MOD}{\text{MOD}}
\newcommand{\PMOD}{\text{PMOD}}

\renewcommand{\epsilon}{\varepsilon}

\title{A Real K3 Automorphism with Most of Its Entropy in the Real Part}
\author{Ethan Cohen}
\address{Department of Mathematics, Yale University}
\email{ethan.cohen@yale.edu}
\begin{document}
	\begin{abstract}
		This article describes the first example of a real K3 surface $X$ admitting a real automorphism $f$ satisfying $h_{top}(f, X(\mathbb{C})) < 2 h_{top}(f, X(\mathbb{R}))$. The example presented is a $(2,2,2)$-surface in $\mathbb{P}^1 \times \mathbb{P}^1 \times \mathbb{P}^1$ defined by the vanishing set of $(1 + x^2)(1 + y^2)(1 + z^2) + 10xyz - 2$, first described by McMullen. Along the way, we develop an ad hoc shadowing lemma for $C^2$ (real) surface diffeomorphisms, and apply it to estimate the location of a periodic point in $X(\mathbb{R})$. This result uses the GNU MPFR arbitrary-precision arithmetic library in C and the Flipper computer program. 
	\end{abstract}
	\maketitle
	
	\tableofcontents
	
		\section{Introduction}
	Let $X$ be a K3 surface equipped with a real structure, and let $f:X\to X$ be a \emph{real automorphism} of $X$; that is, a holomorphic diffeomorphism which commutes with the given antiholomorphic involution $\sigma$ of $X$.
	The fixed point set of $\sigma$, denoted $X(\R)$, is either empty or a (possibly disconnected) compact, real-analytic, totally real submanifold with real dimension two. 
	Trivially, $h_{top}(f, X(\R)) \leq  h_{top}(f, X(\C))$. However, in his ICM survey (\cite{CantatICM}), Cantat noted that there were no known examples for which $h_{top}(f, X(\C)) < 2 h_{top}(f, X(\R))$. The purpose of this paper is to provide such an example. 
	
			\begin{figure}[h]
		\centering
		\includegraphics[width=0.5\textwidth]{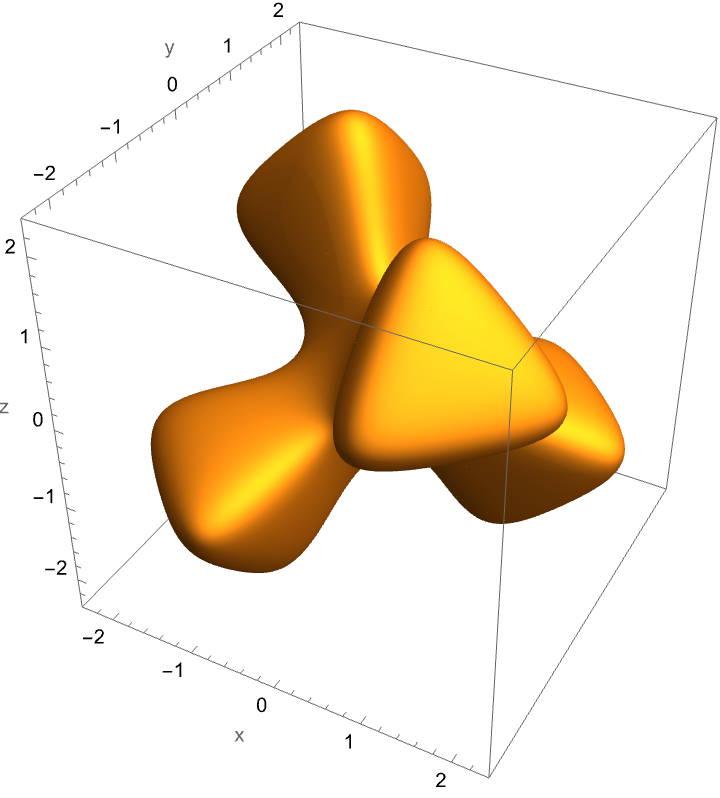}
		\caption{$X_{10}(\R)$, plotted in Mathematica}
		\label{fig:1}
	\end{figure}

	Before beginning, we recall a few related facts. Cantat proved in \cite{MR1864630} that if $Y$ is a connected compact complex surface equipped with an automorphism $f$ of 
	positive topological entropy, then $Y$ is K\"ahler and one of the following holds:
	\begin{enumerate}
		\item  \ either $Y$ is a non-minimal rational surface;
		\item \ or $Y$ is a K3 surface, a complex torus, or an Enriques surface;
		\item \ or there exists a surface $\tilde Y$ from case (2) and a holomorphic birational map $\pi:Y\to \tilde Y$ such that $\tilde f\coloneqq \pi\circ f\circ \pi^{-1}$ is an automorphism of $\tilde Y$.
	\end{enumerate}
	Given a surface $Y$ with a real structure and a real automorphism $f$ with positive entropy, we let $\rho(f)~\coloneqq ~\frac{h_{top}(f, Y(\R))}{h_{top}(f, Y(\C))}$ denote the ratio of the real and complex entropies. 
	
	In case (3), the surface $Y$ is obtained from $\tilde Y$ by a finite sequence of $\tilde f$-equivariant blow-ups along periodic orbits of $\tilde f$. Since the entropy of an automorphism of a K\"ahler manifold is determined by the spectral radius of its action on cohomology (see the Gromov--Yomdin theorem \cite{MR880035}, \cite{MR2026895}, \cite{MR889979}), $h_{top}(f, Y(\C)) = h_{top}(\tilde f, \tilde Y(\C))$. Moreover, even if some of the blow-ups are along an orbit in the real part, one can show 
	$h_{top}(f, Y(\R)) = h_{top}(\tilde f, \tilde Y(\R))$. Therefore, it's enough to focus on cases (1) and (2).
	
	In the case that $Y$ is a rational surface, Diller-Kim produced a family of examples for which $\rho(f) = 1$ (\cite{MR4261942}, Theorem A). They also provided examples for which $\frac 1 2 < \rho(f) < 1$ (\cite{MR4261942}, Theorem B). To my knowledge, in the rational case, no examples are known for which $\rho(f) = 0$. 
	
	For K3 surfaces, it remains open whether $\rho(f) =1$ or $\rho(f) = 0$ can occur. 
	On one hand, Moncet proved that for any $\epsilon>0$, there exists a K3 surface admitting a real automorphism $f$ with $\rho(f)<\epsilon$ (\cite{MR2959025}, Theorem 5.7). He also showed that if $Y$ is a complex torus, then any real automorphism $f$ of $Y$ with positive entropy satisfies $\rho(f) = \frac 1 2$ (\cite{MR2959025}, Proposition 4.4). The same holds when $(Y,f)$ is a so-called Kummer example in the sense of \cite{MR4071328}. For instance, a K3 Kummer example  $(Y,f)$ is obtained by taking $Y$ to be the minimal regular model of the quotient of the standard complex torus $\C^2/\Lambda$ by $(x,y)\mapsto (-x,-y)$, and taking $f$ to be induced by a real hyperbolic automorphism of the torus. 
	
		Generally, a few strategies for bounding $\rho(f)$ from below have emerged. While complex entropy can usually be computed explicitly using the Gromov--Yomdin theorem, the challenge becomes providing a lower bound on the real entropy. In \cite{MR4261942}, Diller-Kim bounded the real entropy by understanding how $f$ acts on the homology of the real part. In forthcoming work, Kim provides sharper lower bounds on $\rho(f)$ for some automorphisms by analyzing homotopy growth rates. Alternatively, Moncet proved that $\rho(f)$ is bounded below by an intrinsic property of the surface $Y$ called its ``concordance'' (\cite{MR2959025}). However, concordance is generally tricky to calculate. 
	In the present work, we bound real entropy by computing the stretch factor of the mapping class of $f$, acting on the real surface with punctures along a real periodic orbit.

	\subsection{Setup and Main Theorem} \label{sec:setup}
	Given nonzero $A\in \R$, consider 
	\[q_A(x,y,z) \coloneqq (1+x^2)(1+y^2)(1+z^2) + Axyz - 2.\] Taking $([x_0: x], [y_0: y], [z_0: z])$ to be homogeneous coordinates for $(\P^1)^3$, 
	one has \[q_A(x,y,z) = \tilde q_A([1: x], [1:y], [1:z])\] where
	\begin{align*}
		\tilde q_A([x_0: x], [y_0: y], [z_0: z])
		\coloneqq (x_0^2 + x^2)(y_0^2 + y^2)(z_0^2 + z^2)
		+ A x y z x_0y_0z_0 - 2x_0^2y_0^2z_0^2.
	\end{align*}
	Let $X_A\subset(\P^1)^3$ denote the hypersurface defined by the vanishing set of $\tilde q_A$. Then $X_A$ is a smooth real projective surface (see Lemma \ref{lem:smoothXA}). A computation using the adjunction formula reveals the canonical bundle of $X_A$ to be trivial. Moreover, the Lefschetz hyperplane theorem implies $b_1(X_A)=0$. It follows that $X_A$ is a K3 surface.
	
	Each $X_A$ is an example of a so-called $(2,2,2)$ surface in $\P^1 \times \P^1\times \P^1$, since $q_A$ has degree two in each variable. There are three natural involutions of any $(2,2,2)$-surface which we recall here. Let $\pi_i:(\P^1)^3\to \P^1$ denote projection to the $i$th coordinate. For $i\neq j$, the map $\pi_{ij} \coloneqq (\pi_i, \pi_j):(\P^1)^3\to (\P^1)^2$ restricted to $X_A$ has fibers generically containing exactly two points and otherwise containing one point or a rational curve. 
	We define $\sigma_k^A$ to be the involution of $X_A$ which preserves each fiber of $\pi_{ij}$ ($i,j\neq k$) and swaps elements of a generic fiber. Since any birational self-map of a K3 surface extends to an automorphism, $\sigma^A_k$ is an automorphism. The three involutions $\sigma_i$ on any $(2,2,2)$-surface generate a free product $\Z/2\Z \star \Z/2\Z \star \Z/2\Z$ in $Aut(X)$ (see \cite{MR1352278}).

	For each nonzero $A\in \R$, we consider the automorphism of $X_A$ given by
	\[f_A \coloneqq  \sigma_3^A\circ \sigma_2^A\circ \sigma_1^A.\] In Section \ref{sec:ComplexEntropy}, we prove that $h_{top}(f_A, X_A(\C))$ is positive and independent of $A$. 
	
	\begin{rem}
		In \cite{MR2959025}, Moncet uses a type of continuity of topological entropy to prove that \[\lim_{A\to 0}h_{top}(f_A, X_A(\R)) = 0.\] 
		He applies a theorem of Cantat from \cite{MR1864630} to obtain that $h_{top}(f_A, X_A(\C))$ is positive and independent of $A$, implying $\lim_{A\to 0}\rho(f_A) = 0$. Our computation of $h_{top}(f_A, X_A(\C))$ differs from Moncet's; he assumes that the fibers of $\pi_{ij}$ do not contain any rational curves, which is not the case. Regardless, $h_{top}(f_A, X_A(\C))$ remains positive and independent of $A$, so Moncet's conclusion holds.
	\end{rem}
	\begin{rem}
	We note that the family of examples $\set{X_A}$ is non-generic within all $(2,2,2)$-surfaces. Cantat proves in \cite{MR1864630} that a generic $(2,2,2)$-surface has Picard number $3$, while we will see that each $X_A$ has Picard number at least $12$. Moreover, Cantat proves that the automorphism $\sigma_3\circ \sigma_2 \circ \sigma_1$ on a generic $(2,2,2)$-surface has entropy $\ln (9 + 4 \sqrt 5)\approx \ln 17.9443$, while we will find $h_{top}(f_A)\approx \ln 6.1393$. 
\end{rem}

	While the calculations in Section \ref{sec:ComplexEntropy} apply for any nonzero $A$, 
	in Section \ref{sec:RealEntropy} we focus on the specific case $A = 10$. The real surface $X_{10}(\R)$ is homeomorphic to a sphere and completely contained in the affine part of $X_{10}$ (see Lemma \ref{lem:star} and Figure \ref{fig:1}).
	Define
	\begin{align}\label{definedQ}
		Q(t)
			\coloneqq t^6-5t^5-6t^4-5t^3-6t^2-5t+1,
	\end{align}
and 
\begin{equation}\label{definedP}
		\begin{aligned}
		P(t) \coloneqq &t^{14}-16t^{13}+78t^{12}-108t^{11}-99t^{10}
		+334t^9+56t^8 \\
		&-556t^7 +56t^6+334t^5-99t^4-108t^3+78t^2-16t+1.
	\end{aligned}
\end{equation}
	Let $\mu$ be the unique real root of $Q$ greater than $1$, and let $\lambda$ be the unique real root of $P$ on the interval $(8.19981448, 8.19981449)$.
	\begin{thm}[Main Theorem]\label{mainTheorem}
		For $X_A, f_A, \lambda,$ and $\mu$ defined above:
		\begin{enumerate}
			\item \ $h_{top}(f_A, X_A(\C))=\ln \mu$ for every nonzero $A\in \R$;
			\item \ there exists a finite $f_{10}$-invariant set $S\subset X_{10}(\R)$ such that the mapping class of $f_{10}^2$ on $X_{10}(\R)\setminus S$ is pseudo-Anosov with stretching factor $\lambda$.
		\end{enumerate}
		As a consequence of (1) and (2), 
		\[	\rho(f_{10}) \geq \frac{\ln\lambda}{2\ln\mu} > \frac 1 2.\]
	\end{thm}
	
	\subsection{Outline}
		Our computation of $h_{top}(f_A,X_A(\C))$ in Section \ref{sec:ComplexEntropy} relies on the Gromov--Yomdin theorem and working with explicit algebraic curves. In Section \ref{sec:ShadowingLemma}, we prove a checkable shadowing lemma for an arbitrary $C^2$-diffeomorphism of a (real) surface. 
	In Sections \ref{sec:PeriodicPoint} and \ref{sec:ApplyShadowing}, we apply the shadowing lemma to a pseudo-orbit of $f_{10}$ on $X_{10}(\R)$, with the help of a computer, to obtain a point $x\in X_{10}(\R)$ satisfying $f_{10}^{10}(x)=x$. Section \ref{sec:arc-certification} certifies that $x,f_{10}(x),\dots,f_{10}^9(x)$ are pairwise distinct; hence $x$ has least period ten. 
	
	In Section \ref{sec:mclass}, we compute the mapping class of $f_{10}^2$ on $X_{10}(\R)\setminus \set{f^i(x) \mid i\in \Z}$ in terms of Dehn half-twists. The mapping class representation of $f_{10}^2$ is written explicitly in Section \ref{section:application}. From there, Flipper (\cite{flipper}) is able to compute the stretch factor (dilatation) of $f_{10}^2$, yielding a lower bound on $h_{top}(f_{10}, X_{10}(\mathbb{R}))$ (see \cite{MR568308}).
	
	\subsection{Accompanying computer code} 
	The GitHub repository accompanying this article is
	\href{https://github.com/ethanhcoo/K3entropy.git}{K3entropy}.
	See its \texttt{README} file for more details. All code is based on	Curt McMullen's `Orbits of Automorphisms of K3 Surfaces' (see \cite{mcmullen2006k3software}). Throughout, when we say that a program \emph{certifies} a statement, we mean that it proves the statement using exact rational arithmetic or outward-rounded interval arithmetic. The outward-rounding procedure, implemented using the GNU MPFR library (see \cite{MR2326955}), is described in Appendix~\ref{sec:errors}. Aside from \texttt{Flipper}, the three programs essential to the main theorem are:
	\begin{itemize}
		\item \texttt{verifyPeriodicOrbit.c}  -- utilized in Section~\ref{sec:ApplyShadowing} -- certifies an upper bound for $\max_i\norm{f_i^c(0)}_2$ (see Lemma \ref{lem:closing} for notation).
		\item  \texttt{derivativesInterval.c} -- utilized in Sections \ref{sec:1} and \ref{sec:kappa} -- certifies that the matrices $\widetilde L_i$ displayed in Table \ref{table:2} satisfy $\norm{Df_i^c(0)-\widetilde L_i}_\infty<10^{-2}$. It also certifies bounds for $\mathcal{D}$ (see Section \ref{sec:charts} for notation) and its first and second partials along a specific pseudo-orbit (Table \ref{table:3}).
	
    \item \texttt{verifyArcs.c} -- used to produce the top panel of Figure~\ref{fig:actual} -- certifies the relative homotopy classes of the images of several arcs under $f_{10}$.

	\end{itemize}

	In Section \ref{sec:algorithm}, we provide an algorithm for realizing the mapping class of a homeomorphism of an $n$-punctured sphere as a product of Dehn half-twists. The algorithm is implemented in \texttt{mclass.c} and takes as input the arc-data obtained by \texttt{verifyArcs.c}, stored in the file \texttt{arc\_certificate\_data.c}. Section~\ref{section:application} explains how the output of \texttt{mclass.c} is illustrated in Figure~\ref{fig:actual}.

	Finally, we use Flipper to compute the stretch factor of the product of Dehn half-twists representing the mapping class of $f^2_{10}$ on $X_{10}(\R)\setminus \set{f^i(x) \mid i\in \Z}$.
	The file \texttt{FlipperDilatation.ipynb} is a Jupyter Notebook which runs Flipper 0.15.6 in SageMath 10.5 (\cite{sagemath}), with \texttt{realalg} 0.3.6. The notebook verifies that the mapping class is pseudo-Anosov, constructs the $42\times42$ hitting matrix for its invariant train track, verifies that $P$ is the minimal polynomial of its stretching factor, and provides an interval containing the stretching factor.

	\subsection{Notation} \label{sec:notations}
	Here we collect some notation used in the rest of the paper:
	\begin{itemize}
		\item[--] For $v=(v_1,\dots,v_d)\in\R^d$, let $\norm{v}_2\coloneqq\left(\sum_{j=1}^d|v_j|^2\right)^{\frac12}$.
		\item[--] For $v=(v_1,\dots,v_d)\in\R^d$, let $\norm{v}_\infty\coloneqq\max_{1\leq j\leq d}|v_j|$.
		We use the same notation for any matrix to denote its maximum absolute entry.
		\item[--] $\norm{\cdot}_0$ is the norm on $(\R^2)^{n}$ defined by $\norm{(\vec v_1, \dots, \vec v_{n})}_0 \coloneqq \max_{1\leq i \leq n}\norm{\vec v_i}_2$.
		\item[--] For $M\in M_{m\times n}(\R)$, $\norm{M}_{op}$ denotes the operator norm for $M$ as a map $(\R^n, \norm{\cdot }_{2})\to (\R^m, \norm{\cdot }_{2})$.
		\item[--] Given $M\in M_{2n\times 2n}(\R)$, $\norm{M}_{op, 0}$ denotes the operator norm of $M$ as a map $((\R^2)^n, \norm{\cdot}_0)\to ((\R^2)^n, \norm{\cdot}_0)$.
		\item[--] For $M\in M_{m\times n}(\R)$, $\norm{M}_F$ is the Frobenius norm.
			\item[--]  For twice-differentiable $f=(f_1,\dots,f_m):\R^n\to \R^m$, $D^2f$ denotes the collection $(D^2f_1,\dots,D^2f_m)$, where the $D^2f_i$ are Hessian matrices.
			\item[--]  For $f:\R^n\to \R^m$, let $\norm{f}_{C^0(U)} \coloneqq \sup_{x\in U} \norm{f(x)}_{2}$.
		\item[--] For $f\in C^k(\R^n)$, let $\norm{f}_{C^k(U),\infty} \coloneqq \max_{|\beta|\leq k}\sup_{x\in U} |D^\beta f(x)|$. 
		\item[--] For a $C^k$-function $f = (f_1, \dots, f_m): \R^n\to \R^m$, let $\norm{f}_{C^k(U),\infty} \coloneqq \max_i\norm{f_i}_{C^k(U),\infty}$.
		\item[--] For $G:\R^n\to M_{\ell \times k}(\R)$, let $\norm{G}_{U, \star} \coloneqq \sup_{x\in U} \norm{G(x)}_{\star}$, where $\star \in \set{op, F, \infty}$.
	\end{itemize}
	
	\subsection{Acknowledgements} 
	The author is grateful to his advisor, Sebastian Hurtado, for his ample guidance and generosity.
	He would also like to thank Serge Cantat for an encouraging conference chat and for his comments on an earlier version of this work.

	\section{Computing $h_{top}(f,X(\C))$} \label{sec:ComplexEntropy}
	\begin{notn}
		For the rest of this section, we fix nonzero $A\in \R$ and let $X\coloneqq X_A$ and $f \coloneqq f_A$.
	\end{notn}
	\subsection{Background}
	The Gromov--Yomdin theorem implies that $h_{top}(f, X(\C)) = \ln R(f^*)$, where $R(f^*)$ denotes the spectral radius of $f^*\in GL(H^*(X, \C))$.
	In this section, we recall a few relevant facts about K3 surfaces (see, for example, \cite{MR3586372}). By Hodge theory, $X$ satisfies 
	\[ \dim H^k(X,\C) = \begin{cases} 
		1& k = 0, 4 \\
		22 & k = 2 \\
		0 & \text{otherwise}.
	\end{cases}
	\]
	Therefore, $f^*$ achieves its spectral radius on $H^2(X,\C)$. 
	Additionally, there exists an $Aut(X)$-invariant decomposition of $H^2(X,\C)$ into its Dolbeault cohomology classes, denoted
	\[H^2(X, \C) = H^{2,0}(X, \C) \oplus H^{1,1}(X, \C) \oplus H^{0,2}(X, \C).\]
	Complex conjugation preserves $H^{1,1}(X, \C)$ and $H^{2,0}(X, \C) \oplus H^{0,2}(X, \C)$, so both are spanned by elements of $H^2(X,\R)$. 
	By the Hodge Index Theorem, the intersection form $\angle{\alpha, \beta} \coloneqq \int_{X}\alpha \wedge  \beta$ is positive-definite on $(H^{2,0}(X, \C) \oplus H^{0,2}(X,\C))\cap H^2(X, \R)$ and has signature $(1, 19)$ on $H^{1,1}(X, \R)\coloneqq H^{1,1}(X, \C)\cap H^2(X,\R)$. 
	Consequently, $f^*$ achieves its spectral radius on $H^{1,1}(X, \R)$. 
	Define the N\'eron-Severi group by
	\[NS(X)\coloneqq H^2(X,~\Z)\cap H^{1,1}(X,~\R)\]
	in the sense that we consider the image of $ H^2(X,~\Z)$ in $ H^2(X,~\R)$ and intersect with $H^{1,1}(X,~\R)$. 
	Then $NS(X)$ is a discrete subgroup of $H^{2}(X, \R)$. Since $Aut(X)$ preserves $H^2(X, \Z)$ and the Hodge decomposition, it preserves $NS(X)$. Since $X$ is a K3 surface, there exist $Aut(X)$-equivariant isomorphisms
		\[D(X)\cong Pic(X) \cong NS(X),\]
	where $Pic(X)$ is the Picard group and
	$D(X)$ is the free abelian group of algebraic curves considered up to linear equivalence. Moreover, under these isomorphisms, $\angle{\cdot,\cdot}$ is carried to the intrinsic intersection forms on $D(X)$ and $Pic(X)$. Since $X$ is projective, the Picard group contains a line bundle of positive self-intersection. Therefore, $\angle{\cdot,\cdot}$ has signature $(1,\rho(X)-1)$ on each of the real vector spaces
	\begin{align*}
		&D(X)\otimes_{\Z}\R,\hspace{1cm}
		Pic(X)\otimes_{\Z}\R,\hspace{1cm}
		NS(X)\otimes_{\Z}\R,
	\end{align*}
	where $\rho(X)\coloneqq \operatorname{rank}Pic(X) =\operatorname{rank}NS(X)$ is the Picard number.
	
	\subsection{Strategy}
	We will find an explicit $\sigma_i$-invariant subspace $W$ of $D(X)\otimes_{\Z}\R$ on which $\angle{\cdot,\cdot}$ is Minkowski. 
	Then $R(f^*) = R(\restr{f^*}{W})$ will be directly computed. 
	
	\subsection{Computations on homology}\label{sec:homology}
	To that end, we consider several algebraic curves in $X$. In homogeneous coordinates $([x_0: x], [y_0: y], [z_0: z])$ for $\P^1\times \P^1\times \P^1$, let
	\begin{align*}
		p_1 &\coloneqq \P^1\times \{[0: 1]\} \times \{[1: i]\}, \ \ 
		p_2 \coloneqq \P^1\times \{[0: 1]\} \times \{[1: -i]\},  \\
		p_3 &\coloneqq \P^1 \times \set{[1: i]} \times \set{[0: 1]}, \ \ 
		p_4 \coloneqq \P^1 \times \set{[1: -i]} \times \set{[0: 1]}, \\
		p_5 &\coloneqq \set{[0:1]} \times \P^1 \times \set{[1:i]},  \ \ 
		p_6 \coloneqq  \set{[0:1]} \times \P^1 \times \set{[1:-i]}, \\
		p_7 &\coloneqq \set{[1:i]} \times \P^1 \times \set{[0:1]}, \ \
		p_8 \coloneqq \set{[1:-i]} \times \P^1 \times \set{[0:1]}, \\
		p_9 &\coloneqq \set{[0: 1]} \times \set{[1:i]} \times \P^1, \ \
		p_{10} \coloneqq  \set{[0: 1]} \times \set{[1:-i]} \times \P^1, \\
		p_{11} &\coloneqq  \set{[1: i]} \times \set{[0:1]} \times \P^1, \ \
		p_{12} \coloneqq  \set{[1: -i]} \times \set{[0:1]} \times \P^1.
	\end{align*}
		Each $p_i$ is birationally equivalent to $\P^1$, and thus has self-intersection $-2$ by the adjunction formula.
	Moreover, each $p_i$ is the unique representative of $[p_i]\in D(X)$, so we use $p_i$ and $[p_i]$ interchangeably.
	Let
	\begin{align*}
		W\coloneqq\operatorname{span}_{\R}\set{p_1,\dots,p_{12}}\subset D(X)\otimes_{\Z}\R.
	\end{align*}
	A calculation, also performed in \cite{Rowe}, yields
	\begin{align*}
		M &\coloneqq (\langle{p_i, p_j}\rangle)_i^j = \begin{pmatrix}
			-2 & 0 & 0 & 0 & 1 & 0 & 0 & 0 & 0 & 0 & 1 & 1 \\
			0 & -2 & 0 & 0 & 0 & 1 & 0 & 0 & 0 & 0 & 1 & 1 \\
			0 & 0 & -2 & 0 & 0 & 0 & 1 & 1 & 1 & 0 & 0 & 0 \\
			0 & 0 & 0 & -2 & 0 & 0 & 1 & 1 & 0 & 1 & 0 & 0 \\
			1 & 0 & 0 & 0 & -2 & 0 & 0 & 0 & 1 & 1 & 0 & 0 \\
			0 & 1 & 0 & 0 & 0 & -2 & 0 & 0 & 1 & 1 & 0 & 0 \\
			0 & 0 & 1 & 1 & 0 & 0 & -2 & 0 & 0 & 0 & 1 & 0 \\
			0 & 0 & 1 & 1 & 0 & 0 & 0 & -2 & 0 & 0 & 0 & 1 \\
			0 & 0 & 1 & 0 & 1 & 1 & 0 & 0 & -2 & 0 & 0 & 0 \\
			0 & 0 & 0 & 1 & 1 & 1 & 0 & 0 & 0 & -2 & 0 & 0 \\
			1 & 1 & 0 & 0 & 0 & 0 & 1 & 0 & 0 & 0 & -2 & 0 \\
			1 & 1 & 0 & 0 & 0 & 0 & 0 & 1 & 0 & 0 & 0 & -2 \\
			\end{pmatrix},
		\end{align*}
	and 
	\begin{align*}
	S\coloneqq (\langle{p_i, \sigma_1 p_j}\rangle)_i^j = \begin{pmatrix}
		-2 & 0 & 0 & 0 & 1 & 0 & 0 & 0 & 0 & 0 & 1 & 1 \\
		0 & -2 & 0 & 0 & 0 & 1 & 0 & 0 & 0 & 0 & 1 & 1 \\
		0 & 0 & -2 & 0 & 0 & 0 & 1 & 1 & 1 & 0 & 0 & 0 \\
		0 & 0 & 0 & -2 & 0 & 0 & 1 & 1 & 0 & 1 & 0 & 0 \\
		1 & 0 & 0 & 0 & 1 & 0 & 0 & 0 & 0 & 0 & 0 & 0 \\
		0 & 1 & 0 & 0 & 0 & 1 & 0 & 0 & 0 & 0 & 0 & 0 \\
		0 & 0 & 1 & 1 & 0 & 0 & 0 & -2 & 0 & 0 & 0 & 1 \\
		0 & 0 & 1 & 1 & 0 & 0 & -2 & 0 & 0 & 0 & 1 & 0 \\
		0 & 0 & 1 & 0 & 0 & 0 & 0 & 0 & 1 & 0 & 0 & 0 \\
		0 & 0 & 0 & 1 & 0 & 0 & 0 & 0 & 0 & 1 & 0 & 0 \\
		1 & 1 & 0 & 0 & 0 & 0 & 0 & 1 & 0 & 0 & 0 & -2 \\
		1 & 1 & 0 & 0 & 0 & 0 & 1 & 0 & 0 & 0 & -2 & 0 \\
	\end{pmatrix}.
\end{align*}
	
		The matrix $M$ is straightforward to compute. For example, the curves $p_1$ and $p_5$ meet transversely at $([0:1],[0:1],[1:i])$, and thus $\angle{p_1, p_5}=1$. To illustrate how to compute $S$, let's calculate $\angle{p_1,\sigma_1(p_5)}$. When $z=[1:i]$, the defining equation becomes $ \tilde q_A([x_0:x],[y_0:y],[1:i])=x_0y_0(A ixy-2x_0y_0).$ For a fixed generic $[y_0:y]$, we have either $x_0 = 0$ or $Aixy=2x_0y_0$. Since $x_0 = 0$ on $p_5$, it follows that $\sigma_1(p_5)=\set{([x_0:x],[y_0:y],[1:i])\mid A ixy=2x_0y_0}.$ This curve meets $p_1$ transversely at $([1:0],[0:1],[1:i])$, so $\angle{p_1,\sigma_1(p_5)}=1.$

	The matrix $M$ is invertible, so $p_1,\dots,p_{12}$ are linearly independent. Diagonalizing $M$ shows that $\langle \cdot,\cdot \rangle$ has signature $(1, 11)$ on $W$.
	Alternatively, we can realize an element of $W$ with positive self-intersection in the following way. Each fiber of $\pi_{i}:X_A\to\P^1$ is a curve given by the equation $\set{x_i = \alpha}$, all representing the same class in $D(X)$, which we denote by $[c_i]$. 
	Notice $\set{x = [0:1]} = p_5  \cup p_6 \cup p_9 \cup p_{10}$ so 
	$[c_1]= p_5  + p_6 + p_9 + p_{10}$. Similarly, $[c_2] = p_1 + p_2 + p_{11} + p_{12}$ and $[c_3] = p_3 + p_4 + p_7 + p_8$. So each $[c_i]$ is contained in $W$. One can also check $\angle{[c_i], [c_j]} = 2$ whenever $i\neq j$ and vanishes otherwise. Thus $[c_i] + [c_j]$, $i\neq j$, has positive self-intersection.

	Define $\tilde \sigma_1:W\to W$ by $\tilde \sigma_1 \coloneqq M^{-1}S$, so that $\tilde \sigma_1$ has the property
	\begin{align}
			\angle{\sigma_1^* w_1, w_2} = \angle{\tilde \sigma_1 w_1, w_2}, \label{eq:angle}
	\end{align} for all $w_1, w_2\in W$. One finds 
	\[\tilde \sigma_1 = \begin{pmatrix}
		1 & 0 & 0 & 0 & -1 & 0 & 0 & 0 & 1 & 1 & 0 & 0 \\
		0 & 1 & 0 & 0 & 0 & -1 & 0 & 0 & 1 & 1 & 0 & 0 \\
		0 & 0 & 1 & 0 & 1 & 1 & 0 & 0 & -1 & 0 & 0 & 0 \\
		0 & 0 & 0 & 1 & 1 & 1 & 0 & 0 & 0 & -1 & 0 & 0 \\
		0 & 0 & 0 & 0 & -1 & 0 & 0 & 0 & 0 & 0 & 0 & 0 \\
		0 & 0 & 0 & 0 & 0 & -1 & 0 & 0 & 0 & 0 & 0 & 0 \\
		0 & 0 & 0 & 0 & 1 & 1 & 0 & 1 & 0 & 0 & 0 & 0 \\
		0 & 0 & 0 & 0 & 1 & 1 & 1 & 0 & 0 & 0 & 0 & 0 \\
		0 & 0 & 0 & 0 & 0 & 0 & 0 & 0 & -1 & 0 & 0 & 0 \\
		0 & 0 & 0 & 0 & 0 & 0 & 0 & 0 & 0 & -1 & 0 & 0 \\
		0 & 0 & 0 & 0 & 0 & 0 & 0 & 0 & 1 & 1 & 0 & 1 \\
		0 & 0 & 0 & 0 & 0 & 0 & 0 & 0 & 1 & 1 & 1 & 0
	\end{pmatrix}.
	\] 
	The fact that $\tilde \sigma_1$ has integer entries and satisfies $\tilde \sigma_1^2 = id$ was not a priori guaranteed and suggests the following claim.
	\begin{claim}
		$\tilde \sigma_1 = \restr{\sigma_1^*}{W}$
	\end{claim}
	\begin{proof}
		Notice that $\tilde \sigma_1^t M \tilde \sigma_1 = SM^{-1}S = M$. 
		Therefore $\angle{\cdot,\cdot}$ is $\tilde\sigma_1$-invariant. 
		Now $\sigma_1^*$ is induced by an automorphism of $X$ and thus also preserves $\angle{\cdot,\cdot}$. 

		Fix $w\in W$ and set $d\coloneqq \sigma_1^*(w)-\tilde \sigma_1(w)$. For every $u\in W$, equation \eqref{eq:angle} gives
		\begin{align*}
			\angle{d,u} =\angle{\sigma_1^*(w),u}-\angle{\tilde \sigma_1(w),u} =0.
		\end{align*}
		Thus $d\in W^\perp$. Moreover, applying \eqref{eq:angle} with $u=\tilde \sigma_1(w)$ and using the invariance of the intersection form gives
		\begin{align*}
			\angle{d,d}
			&=\angle{\sigma_1^*(w),\sigma_1^*(w)}
			+\angle{\tilde \sigma_1(w),\tilde \sigma_1(w)}
			-2\angle{\sigma_1^*(w),\tilde \sigma_1(w)}\\
			&=\angle{w,w}+\angle{w,w}
			-2\angle{\tilde \sigma_1(w),\tilde \sigma_1(w)} =0.
		\end{align*}
		Since $\angle{\cdot,\cdot}$ is negative definite on $W^\perp$, it follows that $d=0$. Therefore $\sigma_1^*(w)=\tilde \sigma_1(w)$ for every $w\in W$, proving the claim.
	\end{proof}
	Let's check that $\sigma_1^*$ does what we expect to $[c_1], [c_2], [c_3]$. 
	Well, $\sigma_1^*$ should fix $[c_2] = p_1 + p_2 + p_{11} + p_{12}$ and $[c_3] = p_3 + p_4 + p_7 + p_8$; our matrix does this.	 Moreover,
	\begin{equation}\label{eqn}
		\pi_{23}^{-1}(\pi_{23}(c_1)) = \sigma_1(c_1)\cup c_1 \cup p_1 \cup p_2\cup p_3\cup p_4. 
	\end{equation}
		Now $\pi_{23}(c_1)\subset \P^{1}\times \P^1$ has bidegree $(2,2)$ and therefore represents the class $2[\set{\alpha} \times \P^1] + 2[\P^1\times \set{\alpha} ]$. Thus \eqref{eqn} becomes  
	\begin{align*}
		2[c_2] + 2[c_3] = \sigma_1^*[c_1]+ [c_1] + p_1 + p_2 + p_3 + p_4,
	\end{align*} and after expansion
	\begin{align*}
			\sigma_1^*(p_5  + p_6 + p_9 + p_{10}) &= 2( p_{11} + p_{12} + p_7 + p_8)  + p_1 + p_2 + p_3 + p_4 - p_5  - p_6 - p_9 - p_{10},
	\end{align*} which our matrix also satisfies. A similar computation is performed in \cite{MR3748233}.

	Let $P_{xy}$ denote the permutation matrix corresponding to
	\begin{align*}
		(1,5)(2,6)(3,7)(4,8)(9,11)(10,12),
	\end{align*}
	and let $P_{xz}$ denote the permutation matrix corresponding to
	\begin{align*}
		(1,11)(2,12)(3,9)(4,10)(5,7)(6,8).
	\end{align*}
	Using symmetry, one obtains
	\[\restr{\sigma_2^*}{W} = P_{xy}\tilde \sigma_1P_{xy}^{-1} = \begin{pmatrix}
		-1 & 0 & 0 & 0 & 0 & 0 & 0 & 0 & 0 & 0 & 0 & 0 \\
		0 & -1 & 0 & 0 & 0 & 0 & 0 & 0 & 0 & 0 & 0 & 0 \\
		1 & 1 & 0 & 1 & 0 & 0 & 0 & 0 & 0 & 0 & 0 & 0 \\
		1 & 1 & 1 & 0 & 0 & 0 & 0 & 0 & 0 & 0 & 0 & 0 \\
		-1 & 0 & 0 & 0 & 1 & 0 & 0 & 0 & 0 & 0 & 1 & 1 \\
		0 & -1 & 0 & 0 & 0 & 1 & 0 & 0 & 0 & 0 & 1 & 1 \\
		1 & 1 & 0 & 0 & 0 & 0 & 1 & 0 & 0 & 0 & -1 & 0 \\
		1 & 1 & 0 & 0 & 0 & 0 & 0 & 1 & 0 & 0 & 0 & -1 \\
		0 & 0 & 0 & 0 & 0 & 0 & 0 & 0 & 0 & 1 & 1 & 1 \\
		0 & 0 & 0 & 0 & 0 & 0 & 0 & 0 & 1 & 0 & 1 & 1 \\
		0 & 0 & 0 & 0 & 0 & 0 & 0 & 0 & 0 & 0 & -1 & 0 \\
		0 & 0 & 0 & 0 & 0 & 0 & 0 & 0 & 0 & 0 & 0 & -1
	\end{pmatrix}\] and
	
	\[\restr{\sigma_3^*}{W} = P_{xz}\tilde \sigma_1P_{xz}^{-1} = \begin{pmatrix}
		0 & 1 & 1 & 1 & 0 & 0 & 0 & 0 & 0 & 0 & 0 & 0 \\
		1 & 0 & 1 & 1 & 0 & 0 & 0 & 0 & 0 & 0 & 0 & 0 \\
		0 & 0 & -1 & 0 & 0 & 0 & 0 & 0 & 0 & 0 & 0 & 0 \\
		0 & 0 & 0 & -1 & 0 & 0 & 0 & 0 & 0 & 0 & 0 & 0 \\
		0 & 0 & 0 & 0 & 0 & 1 & 1 & 1 & 0 & 0 & 0 & 0 \\
		0 & 0 & 0 & 0 & 1 & 0 & 1 & 1 & 0 & 0 & 0 & 0 \\
		0 & 0 & 0 & 0 & 0 & 0 & -1 & 0 & 0 & 0 & 0 & 0 \\
		0 & 0 & 0 & 0 & 0 & 0 & 0 & -1 & 0 & 0 & 0 & 0 \\
		0 & 0 & -1 & 0 & 0 & 0 & 1 & 1 & 1 & 0 & 0 & 0 \\
		0 & 0 & 0 & -1 & 0 & 0 & 1 & 1 & 0 & 1 & 0 & 0 \\
		0 & 0 & 1 & 1 & 0 & 0 & -1 & 0 & 0 & 0 & 1 & 0 \\
		0 & 0 & 1 & 1 & 0 & 0 & 0 & -1 & 0 & 0 & 0 & 1
	\end{pmatrix}.\]
	Thus
	\[\restr{f^*}{W}
	=\restr{\sigma_1^*}{W}\restr{\sigma_2^*}{W}\restr{\sigma_3^*}{W}
	= 
	\begin{pmatrix}
		0 & 0 & 1 & 1 & 0 & -1 & 0 & 0 & 1 & 1 & 1 & 1 \\
		0 & 0 & 1 & 1 & -1 & 0 & 0 & 0 & 1 & 1 & 1 & 1 \\
		0 & 0 & 2 & 2 & 1 & 1 & 0 & 0 & 0 & -1 & 1 & 1 \\
		0 & 0 & 2 & 2 & 1 & 1 & 0 & 0 & -1 & 0 & 1 & 1 \\
		0 & 1 & -1 & -1 & 0 & -1 & 0 & 0 & 0 & 0 & -1 & -1 \\
		1 & 0 & -1 & -1 & -1 & 0 & 0 & 0 & 0 & 0 & -1 & -1 \\
		0 & 0 & 3 & 3 & 1 & 1 & 0 & 0 & 0 & 0 & 2 & 1 \\
		0 & 0 & 3 & 3 & 1 & 1 & 0 & 0 & 0 & 0 & 1 & 2 \\
		0 & 0 & -2 & -1 & 0 & 0 & 0 & 0 & 0 & -1 & -1 & -1 \\
		0 & 0 & -1 & -2 & 0 & 0 & 0 & 0 & -1 & 0 & -1 & -1 \\
		0 & 0 & 2 & 2 & 0 & 0 & 0 & 1 & 1 & 1 & 2 & 1 \\
		0 & 0 & 2 & 2 & 0 & 0 & 1 & 0 & 1 & 1 & 1 & 2
	\end{pmatrix}.
	\]
	Since $\angle{\cdot,\cdot}$ has signature $(1,11)$ on $W$ and
	signature $(1,19)$ on $H^{1,1}(X,\R)$, it is negative definite on
	the orthogonal complement $W^\perp\subset H^{1,1}(X,\R)$. Because
	$W$ is $f^*$-invariant and $f^*$ preserves the intersection form,
	$W^\perp$ is also $f^*$-invariant. Thus $\restr{f^*}{W^\perp}$ is an
	isometry of a negative-definite space and has spectral radius $1$. Hence, $R(f^*)=R\left(\restr{f^*}{W}\right).$
		 We find that the characteristic polynomial of $\restr{f^*}{W}$ is 
	\begin{align*}
			r(x) =Q(x)(x^2-x+1)^3,
	\end{align*}
	where $Q$ is as defined earlier. The factor $x^2-x+1$ is cyclotomic, while $Q$ is a Salem polynomial. In particular, the spectral radius $R(\restr{f^*}{W})$ is the unique real root of $Q$ of	modulus $>1$, namely $\mu\approx 6.1393$. By the Gromov--Yomdin theorem, 
	\begin{align}\label{ComplexEntropy}
		h_{top}(f, X(\C)) = \ln\mu \approx \ln 6.1393.
	\end{align}
	\begin{rem}
		This calculation demonstrates that $h_{top}(f_A)$ is independent of $A\in \R^*$. 
	\end{rem}
	\section{Bounding $h_{top}(f, X_{10}(\R))$} \label{sec:RealEntropy}
	\begin{notn}\label{notation}
		For the rest of the article, we let $X\coloneqq X_{10}$,  $f\coloneqq f_{10}$, $q\coloneqq q_{10}$, and $\tilde q\coloneqq \tilde q_{10}$.
	\end{notn}
	
	\subsection{Strategy}
	Given a surface $S  = S_{g,n}$, recall that the mapping class group of $S$ is  
	\[\MOD^{\pm} (S) \coloneqq \operatorname{Homeo}(S)/\operatorname{Homeo}_0(S),\]
	noting that we allow for orientation-reversing homeomorphisms. 
	Let $\MOD^+(S)$ denote the index-$2$ subgroup of $\MOD^{\pm}(S)$ consisting of classes coming from orientation-preserving homeomorphisms.
	Given a homeomorphism $h$ of $S$, we let $[h]\in \MOD^{\pm}(S)$ denote its mapping class.

	In the following sections, we will obtain an $f$-periodic point $x\in X(\R)$ of least period ten.
	Then we examine the action of $f$ on $S_{0,10}\cong X(\R)\setminus \set{f^i(x)}_{i=0}^{9}$. 
		Our goal is to represent $[f^2]$ as a product of Dehn half-twists in $\MOD^+(S_{0,10})$. 
		Such a representation allows Flipper to show that $[f^2]$ is a pseudo-Anosov mapping class and compute its stretch factor, denoted $\lambda([f^2])$.
	Since each $\sigma_i$ is orientation-reversing on $X(\R)$, $f$ is orientation-reversing, which is why we consider $f^2$.
		Recall that a pseudo-Anosov homeomorphism $\varphi$ has minimal entropy in its mapping class, and moreover
	$h_{top}(\varphi) = \ln \lambda([\varphi])$ (see \cite{MR568308}).
	Therefore,
	\[h_{top}(f) = \frac 1 2 h_{top}(f^2) \geq \frac 1 2 \ln \lambda([f^2]).\]
	
	Now, we turn to approximating such a periodic point $x$.
	Since $f$ has positive entropy on $X(\C)$, a classical theorem of Katok implies that the number of $n$-periodic points in $X(\C)$ grows at least exponentially in $n$ (see \cite{MR573822}). 
	However, finding explicit periodic points remains a challenge, especially those in $X(\R)$. 
	To that end, we prove an ad hoc shadowing lemma for $C^2$ surface diffeomorphisms, and then use a computer to exhibit a pseudo-orbit $(x_i)$ with sufficient recurrence and hyperbolicity to admit a nearby periodic point.
	
	\subsection{A shadowing lemma} \label{sec:ShadowingLemma}
	The following shadowing lemma can be thought of as an ad hoc version of the Anosov shadowing lemma, where we only care about data along a given pseudo-orbit. The lemma requires some type of hyperbolicity in coordinates along the pseudo-orbit, control of the second derivative in coordinates along the pseudo-orbit, and high recurrence. The benefit is that all constants can be computed explicitly.
	
	\begin{lem} \label{lem:closing}
		Let $M$ be a smooth surface. Let 
		\begin{itemize}
				\item 	$h:M\to M$ be a $C^2$ diffeomorphism
			\item $x_0,\dots, x_{n-1}\in M$
				\item $\set{\phi_i:V_i\to M}$ be a set of $C^2$ charts such that $x_i = \phi_i(0)$.
		\end{itemize}
		Assume $h(x_i)\in \phi_{i+1\bmod n}(V_{i+1\bmod n})$ for all $0\leq i < n$. For each such $i$, set
		\begin{itemize}
			\item $U_i\coloneqq
			\set{v\in V_i\mid
				h(\phi_i(v))\in\phi_{i+1\bmod n}(V_{i+1\bmod n})}$
			\item $h_i^c \coloneqq \phi_{i+1\bmod n}^{-1}\circ h\circ\phi_i
			:U_i\to\R^2$
			\item $	L_i \coloneqq (Dh_i^c)_0.$
		\end{itemize}
		Define $L:\R^{2n}\to\R^{2n}$ by
			\begin{align*}
				L \coloneqq \begin{pmatrix}
					0&0&\cdots&0&L_{n-1}\\
					L_0&0&\cdots&0&0\\
					0&L_1&\ddots&\vdots&\vdots\\
					\vdots&&\ddots&0&0\\
					0&\cdots&0&L_{n-2}&0
				\end{pmatrix}.
			\end{align*}

		Assuming $L - I_{2n}$ is invertible (or equivalently $L_{n-1}L_{n-2}\cdots L_0 - I_2$ is invertible),
		let $\kappa\coloneqq \norm{(L - I_{2n})^{-1}}_{op,0}$ (recalling Notation \ref{sec:notations}).
		If there exists $\delta > 0$ such that
		\begin{enumerate}
			\item \ for all $i$, $B_{12\kappa\delta}(0)\subset U_i$ 
			\item \ for all $i$,
			\[12\kappa\delta < \frac{1}{\kappa} \frac{1}{16 \norm{D^2h_i^c}_{C^0(B_{12\kappa\delta}(0)), \infty}}\]
			\item \ for all $i$, $\norm{h_i^c(0)}_2 < \delta$
		\end{enumerate}
		then there exists an $h$-periodic point $p\in M$ satisfying $h^n(p)=p$ and $h^i(p)\in \phi_i(B_{6\kappa\delta}(0))$. 
	\end{lem}
		Our proof mirrors Katok and Hasselblatt's proof of the Anosov shadowing lemma in \cite{MR1326374}, Theorem 6.4.15. We delay it until Appendix \ref{app:closing}. 
	\subsection{Where to look for pseudo-orbits}
	Let $\iota:(\P^1)^3\to (\P^1)^3$ denote reflection through $\set{x = -z}$. Notice that $\tilde q$ is $\iota$-invariant, so $\iota(X) = X$. Furthermore,  $C\coloneqq ~\set{x = -z}~\cap ~X$ is cut out by an equation of bidegree $(4,2)$ in $\P^1 \times \P^1$. Consequently, $C$ has arithmetic genus $3 = (4-1)(2-1)$ and thus self-intersection $4 = 3 \cdot 2 - 2$ by the adjunction formula for a K3 surface. 
	A simple computation reveals that $\iota \circ f = f^{-1} \circ \iota$. Therefore, any $p \in C \cap f^{n}(C)$ is fixed by $f^{2n}$. 
	
	Recall that there is a unique projective class in the isotropic cone of $\angle{\cdot,\cdot}$ in $H^{1,1}(X, \R)$ that is the attracting fixed point for the action of $f$. 
	Let $v_f^+$ denote a representative of this class lying in the closure of the component of the positive cone containing $[C]$. Similarly, $f$ has a repelling projective class from which we choose a representative $v_f^-$ in the closure of the same component, normalized so that $\angle{v_f^+, v_f^-} = 1$.
	Then 
	\[\lim_{n\to \infty }e^{-n h_{top}(f)} \angle{C, f^n(C)} = \angle{C, v_f^-}\cdot \angle{C, v_f^+}\]
	(see \cite{cantatOnline}, section 2.3).
	Since $[C]$ has positive self-intersection and $v_f^\pm$ lie in the closure of the same positive-cone component as $[C]$, $\angle{C, v_f^-}$ and  $\angle{C, v_f^+}$ are positive.
	Thus, the intersection number $\angle{C,f^n(C)}$ grows like $e^{n h_{top}(f)}$. Since each point of $C\cap f^n(C)$ is fixed by $f^{2n}$, these intersections provide a natural place to search for periodic points. A priori, these intersection points need not be contained in $X(\R)$. However, one can plot $C$ and $f^n(C)$ and see many intersections in $X(\R)$; in the following section, we will input into Lemma \ref{lem:closing} a pseudo-orbit coming from an approximate intersection point of $C$ and $f^5(C)$ in $X(\R)$.
	
	It's worth noting that, if $\operatorname{proj}_W$ denotes the orthogonal projection onto $W$ with respect to the intersection form, one can calculate the intersection of $C$ with each $p_i$ and use the intersection matrix from Section \ref{sec:homology} to arrive at
	$\angle{\operatorname{proj}_W [C], \operatorname{proj}_W [C]} = 4$. But $\angle{[C], [C]} = 4$ and $W^\perp$ is negative-definite, so in fact $[C]\in W$. Therefore, our formula for the action of $f^*$ on $W$
	allows us to explicitly compute $\angle{C, f^n(C)}$ for any given $n$. 
	
	\begin{figure}
		\centering
		\includegraphics[width=.65\textwidth]{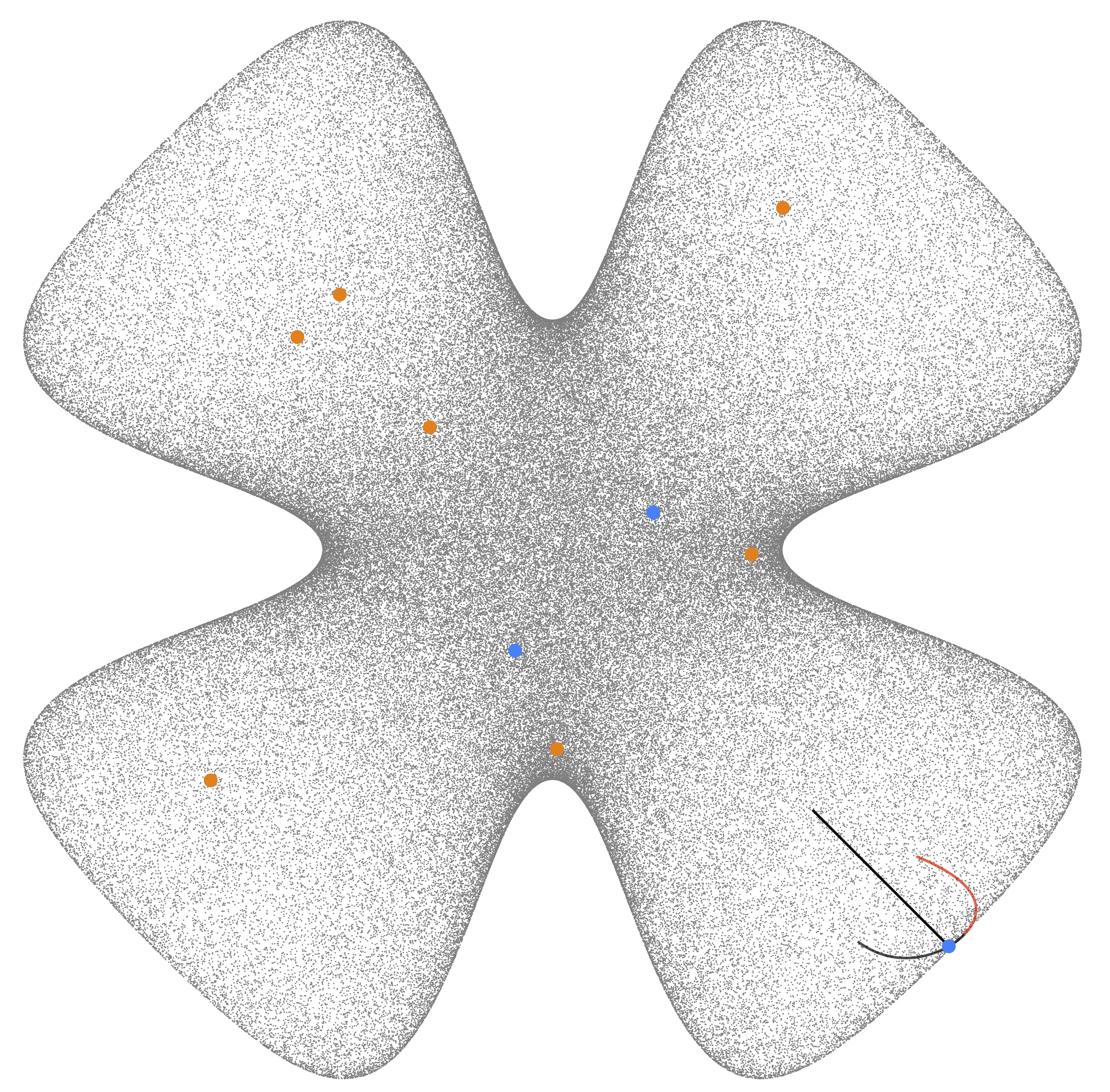} 
			\caption{The pseudo-orbit $x_i = \tilde \phi_i(a_i,b_i)$ from Table \ref{table:1} plotted on $X(\R)$, and projected into the $z$-$x$ plane. Blue dots denote points on the top sheet and orange dots denote points on the bottom sheet. The curve $C \coloneqq \set{x = -z}\cap X$ can be seen intersecting $f^5(C)$ near $x_0$. This figure is obtained using a program developed by C. McMullen.}
		\label{fig:IntersectingCurves}
	\end{figure}
	
	\subsection{A specific pseudo-orbit} \label{sec:PeriodicPoint}
	By the previous section, any intersection between $C$ and $f^n(C)$ in $X(\R)$ has period dividing $2n$. Figure \ref{fig:IntersectingCurves} shows such an intersection when $n = 5$. In the figure, we selected a point near the intersection and plotted the first ten points in its orbit. This will serve as our pseudo-orbit. In the following section, we define the $\phi_i$ used in Lemma \ref{lem:closing}.
	
		\subsubsection{Charts}\label{sec:charts}
			Let $\mathcal{D}(x,y)$ denote the discriminant of $q$ with respect to $z$, i.e.,
		\[\mathcal{D}(x,y) \coloneqq 10^2x^2y^2 + 8 (1+x^2)(1+y^2) - 4(1+x^2)^2(1+y^2)^2.\]
		Let
		\begin{align*}
			\Omega\coloneqq
			\set{(u,v)\in\C^2\mid
			u,v\notin\set{\pm i}
			\text{ and }\mathcal{D}(u,v)\notin(-\infty,0]}.
		\end{align*}
		On $\Omega$, define
		\begin{align*}
			p_{\pm}(u,v)
			\coloneqq
			\frac{-10uv\pm\mathcal{D}(u,v)^{\frac12}}
			{2(1+u^2)(1+v^2)},
		\end{align*}
		where $w\mapsto w^{\frac12}$ denotes the principal branch of the square root, holomorphic on $\C\setminus(-\infty,0]$. Then the maps
		\begin{align*}
			&\Psi_1^{\pm}(y,z) \coloneqq([1:p_\pm(y,z)],[1:y],[1:z]), \hspace{1 cm}
			\Psi_2^{\pm}(x,z)
			\coloneqq([1:x],[1:p_\pm(x,z)],[1:z]),\\
			&\hspace{3.5 cm}
			\Psi_3^{\pm}(x,y) \coloneqq([1:x],[1:y],[1:p_\pm(x,y)])
		\end{align*}
		are inverse holomorphic coordinate charts on $X$ defined on $\Omega$. A direct calculation gives
		\begin{align*}
			(\partial_zq(x,y,z))^2-\mathcal{D}(x,y)
			=4(1+x^2)(1+y^2)q(x,y,z).
		\end{align*}
		Hence, on $X=\set{q=0}$, $\mathcal{D}(x,y)=(\partial_zq(x,y,z))^2.$
		By symmetry, analogous identities hold for the other two partial derivatives of $q$. Therefore, since $X$ is smooth, $\mathcal{D}(x,y)$, $\mathcal{D}(x,z)$ and $\mathcal{D}(y,z)$ cannot vanish simultaneously on $X$. It follows that the restrictions of the three coordinate projections to the six open subsets $\Psi_j^\pm(\Omega\cap\R^2)\subset X(\R)$, for $j=1,2,3$, whose inverses are the corresponding maps $\Psi_j^\pm$, form an atlas of $X(\R)$.
		
	We consider the collection of points $(a_i, b_i)\in \Q^2$ and charts $\tilde \phi_i$ enumerated in Table \ref{table:1}. For each $i$, define $x_i\coloneqq \tilde \phi_i(a_i, b_i).$
	Moreover, let
		\[\phi_i(\zeta, \gamma) \coloneqq  \tilde \phi_i(\zeta + a_i, \gamma + b_i),\] so that the $\phi_i$'s and $x_i$'s are as in Lemma \ref{lem:closing}.
	Here 
	\[V_i = \set{(x,y) \in \R^2 \mid  \mathcal{D}(x + a_i, y+b_i)> 0}\] and 
	\[U_i\coloneqq\set{v\in V_i\mid f(\phi_i(v))\in\phi_{i+1\bmod{10}}(V_{i+1\bmod{10}})}.\]
	Set
	\[f_i^c\coloneqq\phi_{i+1\bmod{10}}^{-1}\circ f\circ\phi_i:U_i\to\R^2.\]
	
		We emphasize that the \emph{exact} values of the $a_i$ and $b_i$ are displayed in Table \ref{table:1}.
	We will see that this pseudo-orbit has enough recurrence and hyperbolicity to apply Lemma \ref{lem:closing}. One could choose $a_i, b_i$ more carefully to get even more recurrence, with a negligible change in the first and second derivatives of $f_i^c$. 
	
	\subsection{Applying the shadowing lemma} \label{sec:ApplyShadowing}
	Recall the notation used in Lemma \ref{lem:closing}.	Fix $\epsilon' \coloneqq 10^{-18}$.
	In Section \ref{sec:1}, we prove that $B_{\epsilon'}(0) \subset U_i$
	and 
		$\max_i\norm{D^2f^c_i}_{C^0(B_{\epsilon'}(0)), \infty} < 1.8\cdot 10^{14}$. Moreover, in Section \ref{sec:kappa} we prove that $\kappa < 83$ for our choices of $\phi_i$.
		In 	\href{https://github.com/ethanhcoo/K3entropy.git}{K3entropy}, \texttt{verifyPeriodicOrbit.c} uses outward-rounded interval arithmetic based on the GNU MPFR library to certify that $0\in U_i$ and $\norm{f^c_i(0)}_2 < 10^{-29}$ for all $i$; see Appendix \ref{sec:errors}.
		Set $\delta\coloneqq10^{-29}$. Then
		\begin{align*}
			12\kappa\delta
			&<10^{-26}
			<\epsilon',
		\end{align*}
		and
		\begin{align*}
			12\kappa\delta
			&<10^{-26}
			<\frac{1}{16\cdot83\cdot1.8\cdot10^{14}}
			<\frac{1}{16\kappa\max_i\norm{D^2f_i^c}_{C^0(B_{\epsilon'}(0)),\infty}}.
		\end{align*}
		Therefore, Lemma \ref{lem:closing} gives a point $x\in X(\R)$ such that $f^{10}(x)=x$ and
		\begin{align}\label{eq:final-shadowing-radius}
			\norm{\phi_i^{-1}(f^i(x))}_2
				<\eta,
		\end{align} for each $0\leq i\leq9$, where $\eta\coloneqq6\kappa\delta$. For the computer-assisted certifications in Section~\ref{sec:arc-certification}, set
		\begin{align*}
			\bar\eta\coloneqq5\cdot10^{-27}.
		\end{align*}
			Then $\eta<\bar\eta$, so $B_\eta(0)\subset[-\bar\eta,\bar\eta]^2$.
	\begin{table}
		\centering
		\begin{tabular}{c|ccc}
			$i$ & $\tilde{\phi}_i$ & $a_i$ & $b_i$ \\
			\hline
			$0$ & $\Psi_1^-$ & 1.041643093944314148360673792017 & 1.726895448754858426328854724474 \\
			$1$ & $\Psi_3^-$ & -0.439586738044637984442175311821 & 0.555943953085459715621476373770 \\
			$2$ & $\Psi_2^-$ & 1.111402054756051352317454938205 & -0.926435350008842162121508383319 \\
			$3$ & $\Psi_1^-$ & -0.328869789067645570794396144391 & 0.867950394543647373540816310310 \\
			$4$ & $\Psi_2^-$ & 1.488818954806569814700993326668 & 1.004464450964796276608907444033 \\
			$5$ & $\Psi_2^-$ & 0.533829900932504729554816817729 & -0.533829900932504729554816817729 \\
			$6$ & $\Psi_2^-$ & -1.004464450964796276608907444033 & -1.488818954806569814700993326668 \\
			$7$ & $\Psi_3^+$ & -0.867950394543647373540816310310 & -0.328869789067645570794396144391 \\
			$8$ & $\Psi_2^-$ & 0.926435350008842162121508383319 & -1.111402054756051352317454938205 \\
			$9$ & $\Psi_1^+$ & 0.555943953085459715621476373770 & 0.439586738044637984442175311821 \\
		\end{tabular}
		\vspace{.25 cm}\\
		\caption{A pseudo-orbit}
		\label{table:1}
	\end{table}
	
	\subsubsection{Bounding $\norm{D^2f^c_i}_{C^0(B_{\epsilon}(0)), \infty}$ for small $\epsilon$} \label{sec:1}
	Define
	\begin{align*}
		K\coloneqq \max_i |\mathcal D (a_i, b_i)|, \hspace{1 cm }R\coloneqq \max_i \norm{D\mathcal D(a_i, b_i)}_{\infty}, \hspace{1 cm} M \coloneqq \max_i \norm{D^2\mathcal D(a_i, b_i)}_{\infty}.
			\end{align*} In \href{https://github.com/ethanhcoo/K3entropy.git}{K3entropy}, the program \texttt{derivativesInterval.c} uses outward-rounded interval arithmetic in the GNU MPFR library to enclose $\mathcal D (a_i, b_i)$ and its first and second partials at $(a_i, b_i)$ for each $i$. The interval enclosures certify that the entries displayed in Table \ref{table:3} are the correct roundings of these values to two decimals. It follows from Table \ref{table:3} that $K < 115$, $R < 163$ and $M < 442$. The program \texttt{derivativesInterval.c} also certifies that $\mathcal D > 0$ on $B_{10^{-3}}(a_i, b_i)$ for all $i$; see Appendix \ref{sec:errors}. In other words, $B_{10^{-3}}(0)\subset  V_i$ for all $i$.
	
	For each $i$, let $\xi_i:\R^3\to\R^2$ be the affine map which agrees with $\phi_i^{-1}$ on $\phi_i(V_i)$. In what follows, allowing an abuse of notation, we consider $f_i^c$ to have domain $V_i$ instead of $U_i$ by setting
	\[ f_i^c\coloneqq \xi_{i+1\bmod 10}\circ f\circ\phi_i.
	\]
	\begin{lem} \label{lem:deriv}
		Fix $\epsilon \coloneqq 10^{-5}$. Then $\max_i\norm{Df^c_i}_{C^1(B_{\epsilon}(0)), \infty} < 1.8\cdot 10^{14}$.
		\begin{proof}[Proof of Lemma  \ref{lem:deriv}]
			In Lemma \ref{lem:well} of the Appendix, 
			we bound $\norm{Df}_{C^0(X(\R)), \infty}$ and $\norm{D^2f}_{C^0(X(\R)), \infty}$, where $f$ is considered to be a function from a neighborhood of $X(\R)\subset \R^3$ into $\R^3$.
			The remaining piece is to bound $\norm{Dp}_{C^0(B_{\epsilon}(a_i, b_i)), \infty}$ and 
			$\norm{D^2p}_{C^0(B_{\epsilon}(a_i, b_i)), \infty}$, and apply the chain rule. 
			We proceed via a series of short claims. 
			
			\begin{claim}  
				\label{claim:4} For all $i$,
				\[\max_{|\beta| = 2} \sup_{B_{\epsilon}(a_i, b_i)}|D^\beta \mathcal{D}| \leq \sqrt 2 \cdot 20000\cdot \epsilon + M.\]
			\end{claim}
			
			\begin{claim} \label{claim:5} For all $i$,
				\[ \max_{|\beta| = 1} \sup_{B_{\epsilon}(a_i, b_i)}|D^\beta  \mathcal{D}| \leq \sqrt 2 ( \sqrt 2 \cdot 20000\cdot \epsilon + M) \epsilon + R.\]
			\end{claim}
			
			\begin{claim} \label{claim:6} For all $i$,
				\[\sup_{B_{\epsilon}(a_i, b_i)}| \mathcal{D}| \leq \sqrt 2(\sqrt 2 ( \sqrt 2 \cdot 20000\cdot \epsilon + M) \epsilon + R) \epsilon + K.\]
			\end{claim}
			
			\noindent Claim \ref{claim:4} follows from Lemma \ref{lem:2} of the Appendix and the mean value theorem, with the $\sqrt{2}$ picked up as a dimensionality constant.
			Claim \ref{claim:5} similarly follows from Claim \ref{claim:4}, and Claim \ref{claim:6} from Claim \ref{claim:5}.
			
			\begin{claim} \label{claim:hm} For all $i$, \[ \norm{D^2p}_{C^0(B_{\epsilon}(a_i, b_i)), \infty} \leq 2.5 \cdot 10^4\] and 
				\[\norm{Dp}_{C^0(B_{\epsilon}(a_i, b_i)), \infty} < 122.\]
			\end{claim}
			\noindent Since $\epsilon=10^{-5}$, Claim \ref{claim:4} gives
			\begin{align*}
				\max_i \max_{|\beta| = 2} \sup_{B_{\epsilon}(a_i, b_i)}|D^\beta \mathcal{D}|
				&<\sqrt 2 \cdot 20000\cdot 10^{-5}+442
				<443.
			\end{align*}
			Claim \ref{claim:5} then gives
			\begin{align*}
				\max_i \max_{|\beta| = 1} \sup_{B_{\epsilon}(a_i, b_i)}|D^\beta \mathcal{D}|
				&<\sqrt 2\cdot443\cdot10^{-5}+163
				<164,
			\end{align*}
			so
			\begin{align*}
				\max_i\norm{D\mathcal{D}}_{C^1(B_{\epsilon}(a_i,b_i)),\infty}<443.
			\end{align*}
			Finally, Claim \ref{claim:6} gives
			\begin{align*}
				\max_i\sup_{B_{\epsilon}(a_i,b_i)}|\mathcal{D}|
				&<\sqrt 2\cdot164\cdot10^{-5}+115
				<116.
			\end{align*}
			Moreover, the same mean-value estimate and Table \ref{table:3} give $\mathcal D(x,y)>8$ for every $(x,y)\in B_\epsilon(a_i,b_i)$ and every $i$.
			Then Lemma \ref{lem:7} of the Appendix gives
			\begin{align*}
				\norm{Dp_\pm}_{C^0(B_{\epsilon}(a_i,b_i)),\infty}
				\leq 5+\frac12\sqrt{116}+\frac14\cdot443<122
			\end{align*} and 
		\begin{align*}
			\norm{D^2p_\pm}_{C^0(B_{\epsilon}(a_i,b_i)),\infty}
			\leq10+\sqrt{116}+\frac34\cdot443+\frac18\cdot443^2
			<2.5\cdot10^4,
		\end{align*}
			proving Claim \ref{claim:hm}.
			
			To conclude the proof of Lemma \ref{lem:deriv}, we examine the action of $f^c_i$ on the space of $2$-jets from $\R^2$ to $\R$ (see Appendix \ref{sec:jets}). Indeed, this action is given by a matrix-valued map $J^2(f^c_i):B_\epsilon(0)\to M_{5\times 5}(\R)$ whose entries contain the first and second partials of $f^c_i$. Therefore we have the bound
			\begin{align*}
				\norm{Df^c_i}_{C^1(B_{\epsilon}(0)), \infty} \leq \norm{J^2(f^c_i)}_{C^0(B_{\epsilon}(0)), \infty}.
			\end{align*}
			Jets behave nicely under composition, meaning
			\begin{align*}
				J^2_x(f^c_i)
				=J^2_x(\phi_i)\cdot J^2_{\phi_i(x)}(f)
				\cdot J^2_{f\circ\phi_i(x)}(\xi_{i+1}).
			\end{align*}
			For matrices $A$ and $B$ where $A$ has $r$ columns and $B$ has $r$ rows, the maximum-entry norm satisfies
			\begin{align*}
				\norm{AB}_\infty\leq r\norm{A}_\infty\norm{B}_\infty.
			\end{align*} 
			Now  $J^2_x(\phi_i)$ is $5\times9$, $J^2_{\phi_i(x)}(f)$ is $9\times 9$, and $J^2_{f\circ\phi_i(x)}(\xi_{i+1})$ is $9\times 5$. Hence
			\begin{align*}
				\norm{J^2_x(f^c_i)}_\infty\leq 9\norm{J^2_x(\phi_i)}_\infty\norm{J^2_{\phi_i(x)}(f)
					\cdot J^2_{f\circ\phi_i(x)}(\xi_{i+1})}_\infty.
			\end{align*}
			But $\xi_{i+1}$ is an affine coordinate projection followed by a translation, so $J^2_{f\circ\phi_i(x)}(\xi_{i+1})$ is an inclusion matrix, and 
			\begin{align*}
				\norm{J^2_{\phi_i(x)}(f)
					\cdot J^2_{f\circ\phi_i(x)}(\xi_{i+1})}_\infty \leq \norm{J^2_{\phi_i(x)}(f)}_\infty.
			\end{align*} 
			Moreover, by \eqref{inequality}, $J^2(f)$ satisfies
			\begin{align*}
				\norm{J^2(f)}_{X(\R), \infty}
				\leq
				\max\left(
					\norm{Df}_{C^0(X(\R)), \infty},
					\norm{D^2f}_{C^0(X(\R)), \infty},
					2\norm{Df}^2_{C^0(X(\R)), \infty}
				\right).
			\end{align*} 
			By Lemma \ref{lem:well}, $\norm{Df}_{C^0(X(\R)), \infty}\leq 9\cdot 10^3$ and $\norm{D^2f}_{C^0(X(\R)), \infty} \leq 9^2(2\cdot 10^2)^3$. It follows that
			\begin{align*}
				\norm{J^2(f)}_{X(\R),\infty}
				\leq9^2(2\cdot10^2)^3.
			\end{align*}
			
			Finally, $J^2(\phi_i)$ satisfies
			\begin{align*}
				\norm{J^2(\phi_i)}_{C^0(B_{\epsilon}(0)), \infty}
				\leq
				\max\left(
					1,
					\norm{Dp}_{C^0(B_{\epsilon}(a_i,b_i)), \infty},
					\norm{D^2p}_{C^0(B_{\epsilon}(a_i,b_i)), \infty},
					2\max\left(1,\norm{Dp}_{C^0(B_{\epsilon}(a_i,b_i)), \infty}\right)^2
				\right),
			\end{align*}
			and by Claim \ref{claim:hm}, the right-hand side is bounded above by
			\begin{align*}
				\max(2\cdot122^2,2.5\cdot10^4)<3\cdot10^4.
			\end{align*}
			Consequently,
			\begin{align*}
				\norm{Df^c_i}_{C^1(B_{\epsilon}(0)),\infty}
				\leq\norm{J^2(f^c_i)}_{C^0(B_{\epsilon}(0)),\infty}
				&\leq9\norm{J^2(\phi_i)}_{C^0(B_{\epsilon}(0)),\infty}
				\norm{J^2(f)}_{X(\R),\infty}\\
				&<9\cdot3\cdot10^4\cdot9^2(2\cdot10^2)^3\\
				&<1.8\cdot10^{14},
			\end{align*} proving Lemma  \ref{lem:deriv}.
		\end{proof}
	\end{lem} 
	
	\begin{lem}\label{lem:domain}
		Fix $\epsilon' \coloneqq 10^{-18}$. Assume that $0\in U_i$ for all $i$ and that $\max_i \norm{f_i^c(0)}_2 < \frac 1 2 10^{-3}$. Then  $B_{\epsilon'}(0) \subset U_i$ for all $i$.
		\begin{proof}[Proof of Lemma \ref{lem:domain}]
			By Lemma \ref{lem:deriv}, we have
			\begin{align}
				\max_i\norm{Df^c_i}_{C^0(B_{\epsilon'}(0)), \infty} < 1.8\cdot 10^{14},
			\end{align}
			which implies $\max_i\norm{Df^c_i}_{C^0(B_{\epsilon'}(0)), op} < 2\cdot1.8\cdot 10^{14}$. We chose $\epsilon'$ such that \[\epsilon'\cdot2\cdot1.8\cdot10^{14}<\frac12 10^{-3}.\] Then for each $i$, the image of the $\epsilon'$-ball under $f_i^c$ satisfies
			\[f^c_i(B_{\epsilon'}(0))\subset B_{ \frac 1 2 10^{-3}}(f_i^c(0)).\] By assumption, $ \norm{f_i^c(0)}_2 < \frac 1 2 10^{-3}$, so it follows that \[f^c_i(B_{\epsilon'}(0))\subset B_{ 10^{-3}}(0)\subset V_{i+1\bmod{10}}.\] Now
			\begin{align*}
				\xi_{i+1\bmod{10}}^{-1}
				\left(B_{ \frac 1 2 10^{-3}}(f_i^c(0))\right)\cap X(\R)
			\end{align*}
			has two connected components, one on each of the two sheets over $B_{ \frac 1 2 10^{-3}}(f_i^c(0))$. The connected set $f\circ\phi_i(B_{\epsilon'}(0))$ is contained in their union. Since $0\in U_i$, this set meets the component corresponding to $\phi_{i+1\bmod{10}}$, and hence is contained entirely in that component. Therefore, $B_{\epsilon'}(0)\subset U_i$.
		\end{proof} 
	\end{lem}
	\subsubsection{Bounding $\kappa$} \label{sec:kappa}
		Recall the definitions of $L_i$, $L$ from Lemma \ref{lem:closing}. 
	We hope to bound
	\begin{align*}
		\kappa\coloneqq\norm{(L-I_{20})^{-1}}_{op,0}.
	\end{align*}
		Let $\widetilde L_i\in M_{2\times2}(\Q)$ be the matrices displayed in Table \ref{table:2}.
		 Then the program \texttt{derivativesInterval.c} in \href{https://github.com/ethanhcoo/K3entropy.git}{K3entropy} certifies that
	\begin{align*}
		\norm{L_i-\widetilde L_i}_\infty<10^{-2}.
	\end{align*}
		Let $\widetilde L$ be the cyclic block matrix formed from the $\widetilde L_i$ which is the analog of $L$ in Lemma \ref{lem:closing}. Write $(\widetilde L-I_{20})^{-1}=(B_{ij})_{0\leq i,j<10}$ in $2\times2$ blocks. Since $\norm{B_{ij}}_{op}\leq\norm{B_{ij}}_F$, the triangle inequality applied to each block row gives
	\begin{align*}
	\norm{(\widetilde L-I_{20})^{-1}}_{op,0} \leq\max_{0\leq i<10}\sum_{j=0}^9\norm{B_{ij}}_{op} \leq\max_{0\leq i<10}\sum_{j=0}^9\norm{B_{ij}}_F
		<31.13,
\end{align*}
	where the final inequality is a direct calculation using the exact rational entries of $(\widetilde L-I_{20})^{-1}$.
	
	Each entry of $L_i-\widetilde L_i$ has absolute value less than $10^{-2}$, so
	\begin{align*}
		\norm{L_i-\widetilde L_i}_{op}
		&\leq\norm{L_i-\widetilde L_i}_F
		<2\cdot10^{-2}.
	\end{align*}
	Because $L-\widetilde L$ has exactly one nonzero $2\times2$ block in each block row,
	\begin{align*}
		\norm{L-\widetilde L}_{op,0}<2\cdot10^{-2}.
	\end{align*}
		Since $\norm{(\widetilde L - I_{20})^{-1}}_{op, 0} \norm{L - \widetilde L}_{op, 0} < 31.13\cdot2\cdot10^{-2}<1$, the Neumann-series estimate gives
	\begin{align*}
		\kappa =\norm{(L-I_{20})^{-1}}_{op,0}
		\leq
		\frac{\norm{(\widetilde L-I_{20})^{-1}}_{op,0}}
		{1-\norm{(\widetilde L-I_{20})^{-1}}_{op,0}\norm{L-\widetilde L}_{op,0}}<
		\frac{31.13}{1-31.13\cdot2\cdot10^{-2}} <83.
	\end{align*}
	\subsection{Computing the mapping class of $f^2$ on a sphere with punctures}\label{sec:mclass} In Lemma \ref{lem:star}, we demonstrate that radial projection $R:X(\R)\to S^2$, where $S^2\subset \R^3$ denotes the unit sphere, is a diffeomorphism. In a small abuse of notation, we still use $f$ to denote the associated conjugate diffeomorphism of $S^2$. Next, we stereographically project through $\gamma \coloneqq (0,0,1)$ to the $x$-$y$ plane, which we view topologically as an open disk $\mathring D$. 
	Let
	$\psi_\gamma:S^2\setminus\{\gamma\}\to \R^2 \cong \mathring{D}$
	denote stereographic projection through $\gamma$. Section \ref{sec:arc-certification} certifies that the points $R(f^j(x))$, for $0\leq j\leq9$, are pairwise distinct and avoid $\gamma$. We can therefore consider their images under $\psi_\gamma$ in $\mathring D$. These are ten points, which we label from left to right by $z_0, z_1, \dots, z_9$, writing each as $z_i = (u_i, v_i)$ for $u_i < u_{i+1}$.
	Let $\tau$ be the permutation such that
	\begin{align*}
	    z_{\tau(j)} = (\psi_\gamma\circ R)(f^j(x)).
	\end{align*}
	Thus $\psi_\gamma \circ f\circ \psi_\gamma^{-1}$ maps $z_{\tau(j)}$ to $z_{\tau(j+1\bmod 10)}$. 
	
	For each $0\leq i\leq 8$, let $p_i:[0,1]\to\mathring D$ be the oriented straight-line arc from $z_i$ to $z_{i+1}$, parametrized by $p_i(t) = (1-t)z_i+tz_{i+1}$. Let $\hat z_i$ and $\hat p_i$ denote the lifts of $z_i$ and $p_i$, respectively, to $S^2$. Let $[\hat s_i]\in \MOD^+(S_{0,10})$ denote the clockwise half-twist around $\hat p_i$ (coming from the  orientation induced on $S^2$ by the standard orientation in $\R^3$), and recall that the half-twists $[\hat s_i^{\pm}]$ generate $\MOD^+(S_{0,10})$. 
	Our goal is to realize $[f^2]\in\MOD^+(S_{0,10})$ as a product of the $[\hat s_{i}^{\pm}]$'s. 
	
	Let $[s_i]\in \MOD^+(\mathring{D}_{10})$ denote the counter-clockwise half-twist about $p_i$, and notice that
	$[\psi_\gamma^{-1}\circ s_i\circ \psi_\gamma] = [\hat s_i]$ once we extend $\psi_\gamma^{-1}\circ s_i\circ \psi_\gamma$ to $S_{0,10}$ by fixing $\gamma$. Since $f(\gamma)\neq\gamma$ but $f^2(\gamma)=\gamma$, the map $f$ does not directly descend to a homeomorphism of $\mathring{D}_{10}$. However, there exists $f'\in [f]\in \MOD^\pm(S_{0, 10})$ such that $f'(\gamma) = \gamma$ and $f'\circ \hat p_i = f\circ \hat p_i$ as parametrized paths for all $i$. Indeed, this follows from the fact that the union $\cup_i \hat p_i$ does not separate $S^2$ and does not contain $\gamma$ or $f^{-1}(\gamma)$. We use $f'$ to denote the induced homeomorphism of $\mathring{D}_{10}\cong S_{0, 11}$. 
	Then \[f'\circ p_i= \psi_\gamma\circ f\circ\hat p_i.\]

	Let $S\coloneqq\set{z_0,z_1,\dots,z_9}$.
	\begin{defn}\label{def:arc-data}
		Let $T\subset\R^2$ be a finite set of marked points with distinct first coordinates, and let $\alpha:[0,1]\to\R^2$ be a continuous path whose endpoints lie in $T$ and such that $\alpha(t)\notin T$ for all $t\in(0,1)$. Assume that there are only finitely many $t\in(0,1)$ for which $\alpha(t)$ lies on a vertical line through a point of $T$, and that at each such $t$, the path passes from one side of that vertical line to the other. The \emph{arc-data of $\alpha$ with respect to $T$} is the sequence of marked points whose vertical lines $\alpha$ crosses, with repetitions recorded, together with whether $\alpha$ passes above or below each marked point.
	\end{defn}
	\begin{rem}
		If two paths with the same endpoints have the same arc-data with respect to $T$, then they are homotopic relative to $T$. Indeed, the vertical lines through the points of $T$ divide $\R^2$ into simply connected regions. After sliding corresponding crossing points along the same components of these vertical lines, the portions of the two paths between corresponding crossings can be homotoped to one another within these regions.
	\end{rem}
	Arc-data serves as an intermediate combinatorial tool for determining the relative homotopy classes of the paths $f'(p_i)$. Indeed, in Section \ref{sec:arc-certification} we construct rational polygonal paths in a combinatorially isomorphic marked-point system and certify that their relative homotopy classes correspond to those of the paths $f'(p_i)$. We compute the arc-data of these polygonal paths exactly, perform homotopy-preserving simplifications, and use the resulting sequences as input to \texttt{mclass.c}. The top panel of Figure \ref{fig:actual} depicts the certified relative homotopy classes of the paths $f'(p_i)$.
	
	In what follows, we use $s_i.s_k$ to denote the composition $s_i\circ s_k$\footnotemark.
	\footnotetext{This convention is consistent with the Flipper program.}
	Using these simplified sequences, the algorithm described in Section \ref{sec:algorithm}	finds $i_k\in \set{0,1,\dots, 8}$ and $n_k\in \set{1,-1}$ such that
	\[g\coloneqq s_{i_0}^{n_0}.s_{i_1}^{n_1}. \dots .s_{i_\ell}^{n_\ell}\in \operatorname{Homeo}^+(\mathring{D}_{10})\] maps each $f'(p_i)$ to $p_i$ up to homotopy in $\mathring{D}_{10}$. This algorithm is implemented in	\href{https://github.com/ethanhcoo/K3entropy.git}{K3entropy}, \texttt{mclass.c}. The program is used to compute the words $g_i$ displayed in Figure \ref{fig:actual}, but the proof does not rely on its correctness, since each word can be checked directly from the figure by tracing the indicated half-twists.

	Then 
	\[\hat g\coloneqq \hat s_{i_0}^{n_0}.\hat s_{i_1}^{n_1}. \dots .\hat s_{i_\ell}^{n_\ell}\in \operatorname{Homeo}^+(S_{0,10}\setminus \set{\gamma})\]
	maps each $f(\hat p_i)$ to $\hat p_i$ up to homotopy in $S_{0, 11} = S_{0,10}\setminus \set{\gamma}$.
	Thus $[\hat g\circ f]\in \MOD^{\pm}(S_{0,10})$ is an orientation-reversing mapping class which fixes each $\hat p_i$; there is only one such mapping class.
	Indeed, let $\nu$ be the unique (up to homotopy in $S_{0,10}$) arc from $\hat z_9$ to $\hat z_0$ that's disjoint from each $\hat p_i$, and let
	$L$ be the closed loop defined by the concatenation
	$\hat p_0.\hat p_1. \dots .\hat p_8. \nu$. Let $\mathcal R\in \operatorname{Homeo}(S_{0,10})$ be a reflection through $L$.
	Notice $\hat g\circ f\circ \mathcal R$ is orientation-preserving and preserves the homotopy class of each $\hat p_i$, so its mapping class is trivial.
	Thus $[\hat g\circ f] = [\mathcal R].$
	
	We are not quite done, since we wanted to realize $[f^2]\in\MOD^+(S_{0,10})$ as a product of 
	the $[\hat s_{i}^{\pm}]$'s. Notice that $[\mathcal R\circ \hat s_i^{\pm}] = [\hat s_i^{\mp}\circ \mathcal R]$ for all $i$. Therefore 
	\begin{equation}\label{fSquared}
		\begin{aligned}
			f^2 &\cong \hat g^{-1} \circ \mathcal R \circ \hat g^{-1} \circ \mathcal R \\
			&\cong (\hat{s}_{i_\ell}^{-n_\ell}  . \hat{s}_{i_{\ell-1}}^{-n_{\ell-1}}.  \dots . \hat{s}_{i_0}^{-n_0})\circ (\hat{s}_{i_\ell}^{n_\ell} . \hat{s}_{i_{\ell-1}}^{n_{\ell-1}}. \dots . \hat{s}_{i_0}^{n_0}) \\
			&\cong \hat g^{-1}\circ \overline{g},
		\end{aligned}
	\end{equation}
	where $\overline{g}$ reverses the order of the product decomposition of $\hat g$, meaning
	\begin{align*}
		\overline{g}\coloneqq \hat{s}_{i_\ell}^{n_\ell} . \hat{s}_{i_{\ell-1}}^{n_{\ell-1}}. \dots . \hat{s}_{i_0}^{n_0}.
	\end{align*}

	\subsubsection{An algorithm for computing $g$}\label{sec:algorithm}
	Let $f', z_i, p_i$, and $s_i^\pm$ be as above, but consider the more general case of $\mathring{D}_n$. In what follows, we consider the $z_i$ as marked points so that homeomorphisms of $\mathring{D}_n$ permute the marked points. 
	
	Let $\MOD^+_{0}(\mathring{D}_n)\subset \MOD^+(\mathring{D}_n)$ denote the subgroup generated by homotopy classes of homeomorphisms that fix the marked point $z_0$. The following lemma will be useful later:
	\begin{lem}\label{lem:generators}
		The subgroup $\MOD^+_{0}(\mathring{D}_n)$ satisfies
		$$\MOD^+_{0}(\mathring{D}_n) = \angle{[\tau_{p_0}], [s_1], [s_2], \dots, [s_{n-2}]}$$ where $\tau_{p_0} = s_0.s_0$ is a full twist around $p_0$.  
		\begin{proof}[Proof of Lemma \ref{lem:generators}]
			For each $0\leq i < j \leq n-1$, define 
			\[A_{i,j} \coloneqq (s_{j-1}. s_{j-2}. \dots
			s_{i+1}).s_i^2. (s_{j-1}. s_{j-2}. \dots.  s_{i+1})^{-1}.\]
			Then $A_{i,j}$ represents a Dehn twist around some loop containing just $z_i$ and $z_j$, and notice that $A_{0,1} = \tau_{p_0}$. Moreover, the collection $\set{A_{i,j}\mid 0\leq i < j \leq n-1}$ generates the pure mapping class group $\PMOD^+(\mathring{D}_n)$ (see \cite{MR2435235}). Therefore $\PMOD^+(\mathring{D}_n) \leq \angle{[\tau_{p_0}], [s_1], [s_2], \dots, [s_{n-2}]}$. Finally, any  $[\phi]\in \MOD_0^+(\mathring{D}_n)$ can be written as a product of $[\phi_1]\in \PMOD^+(\mathring{D}_n)$ and $[\phi_2]\in \angle{[s_1], [s_2], \dots, [s_{n-2}]}$, proving the lemma.
		\end{proof}
	\end{lem}

\begin{figure}
		\includegraphics[width=1\textwidth]{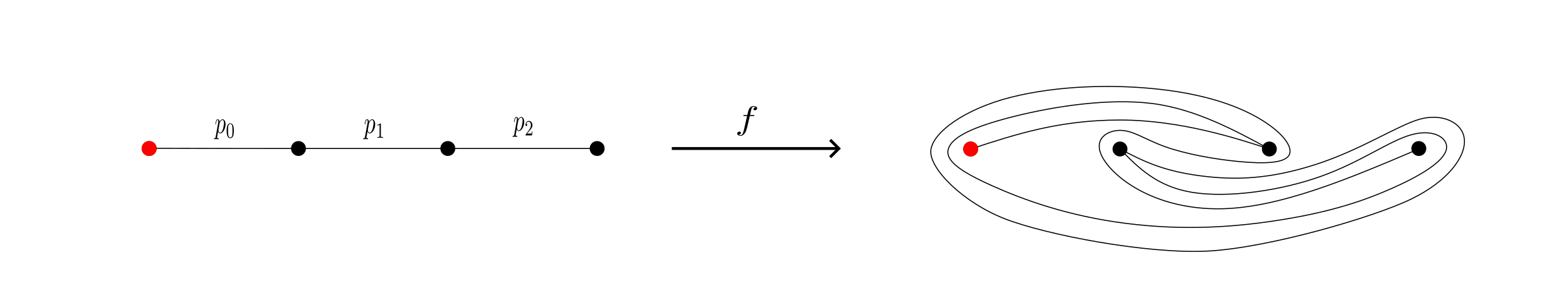}
		\vspace{-.8 cm} \\
		\rule{\linewidth}{0.5pt} 
		\vspace{-.7 cm}
		\begin{multicols}{3}
			$i = 0$\par 
			$i = 1$\par 
			$i = 2$\par 
		\end{multicols}
		\vspace{-.5 cm} 
		\rule{\linewidth}{0.5pt} 
		\vspace{-.8 cm} 
		\begin{multicols}{3}
			\includegraphics[width=\linewidth]{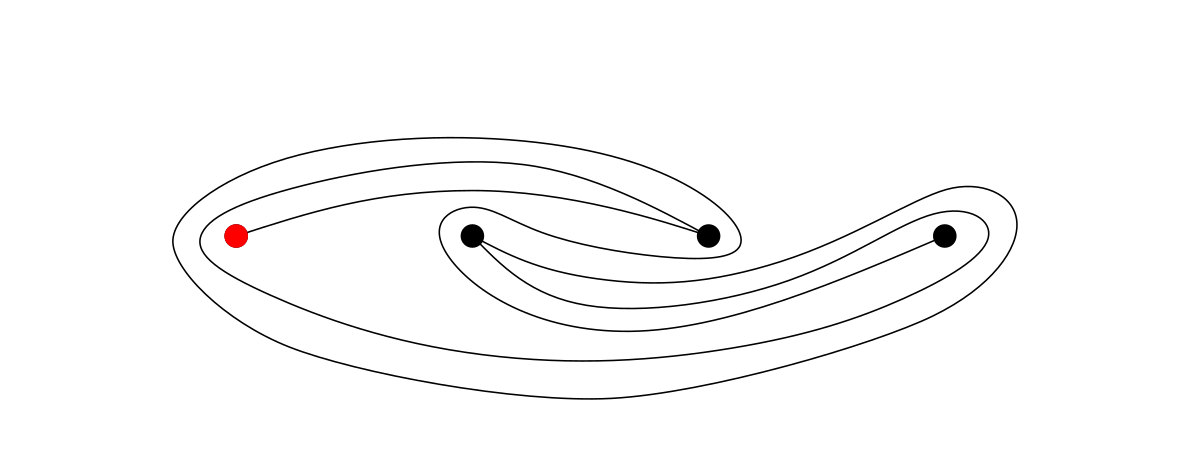}\par 
			\includegraphics[width=\linewidth]{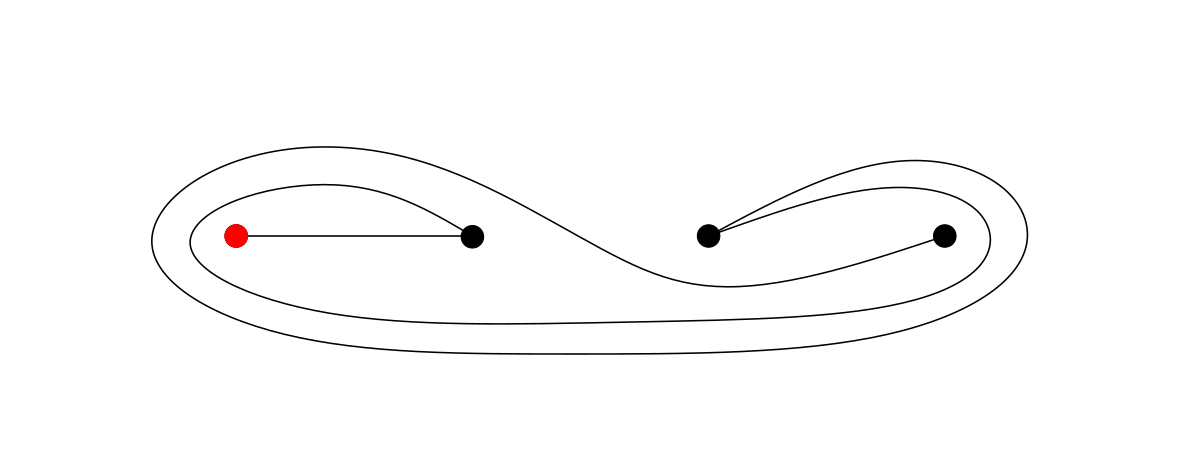}\par 
			\includegraphics[width=\linewidth]{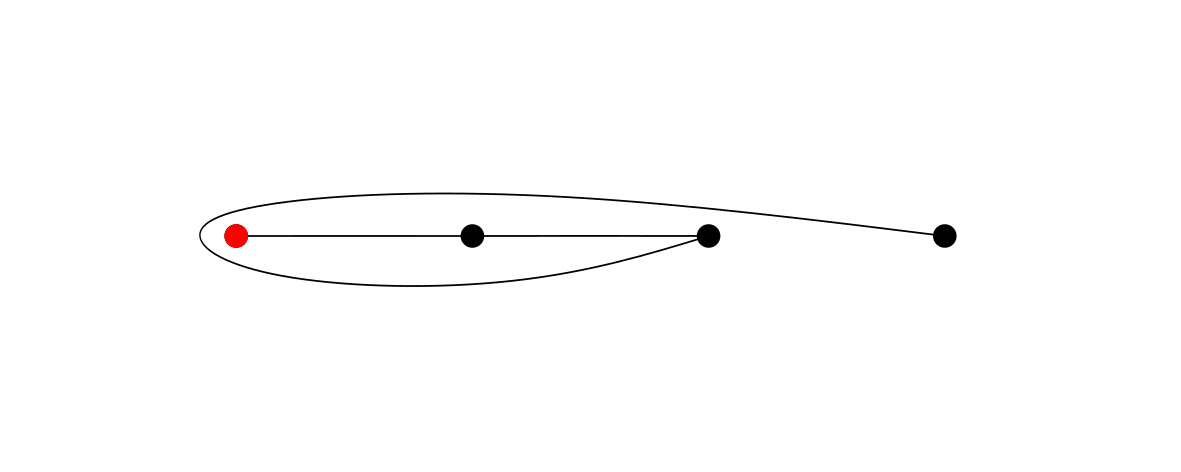}\par 
		\end{multicols}
		\vspace{-1 cm} 
		\begin{multicols}{3}
			\makebox[\textwidth][l]{\textbf{base case:} \hspace{0 cm} $g_0 = s_1$}\par 
			\makebox[\textwidth][l]{\textbf{step 1:} \hspace{0 cm} \text{contract $p_0$}}\par 
			\makebox[\textwidth][l]{\textbf{step 1:} \hspace{0 cm} \text{contract $p_0.p_1$}}\par 
		\end{multicols}
		\vspace{-1 cm} 
		\begin{multicols}{3}
			\includegraphics[width=\linewidth]{ex3.png}\par 
			\includegraphics[width=\linewidth]{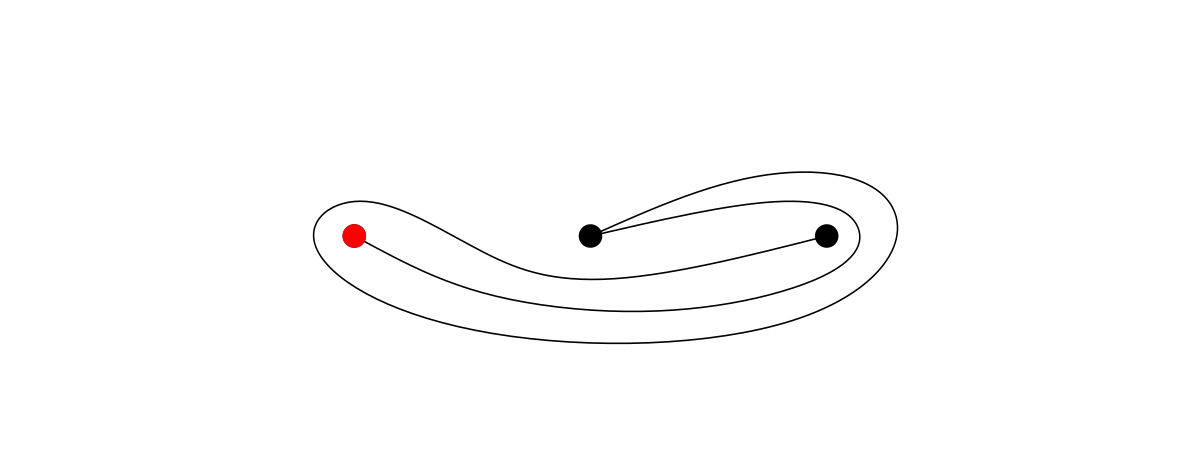}\par 
			\includegraphics[width=\linewidth]{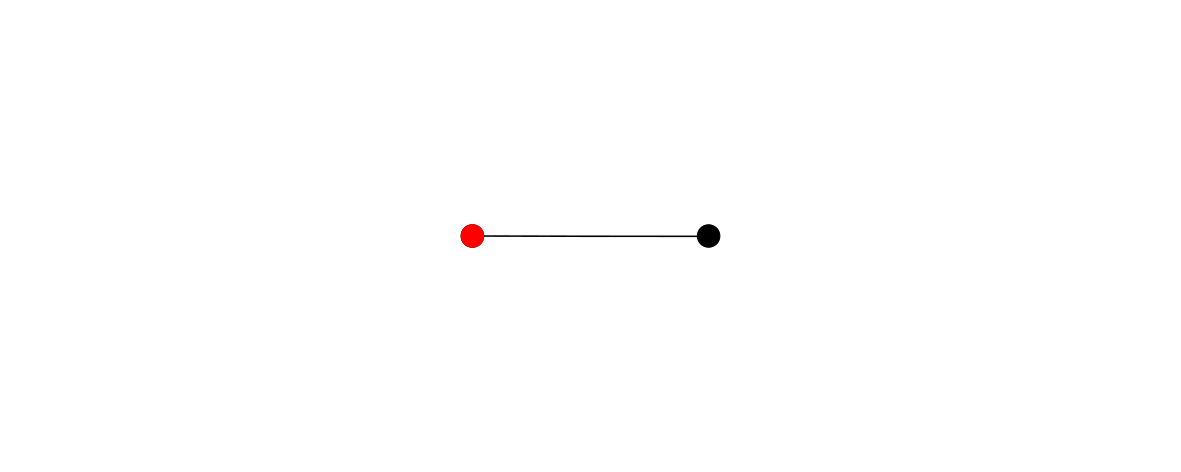}\par 
		\end{multicols}
		\vspace{-1 cm} 
		\begin{multicols}{3}
			\hspace{.1 cm}\par 
			\makebox[\textwidth][l]{\textbf{step 2:} \hspace{0 cm} $\tilde g_1 = s_2^{-1}.s_2^{-1}$}\par 
			\makebox[\textwidth][l]{\textbf{step 2:} \hspace{0 cm} \text{$\tilde g_2 = id$}}\par 
		\end{multicols}
		\vspace{-1 cm} 
		\begin{multicols}{3}
			\hspace{.1 cm} \par   
			\includegraphics[width=\linewidth]{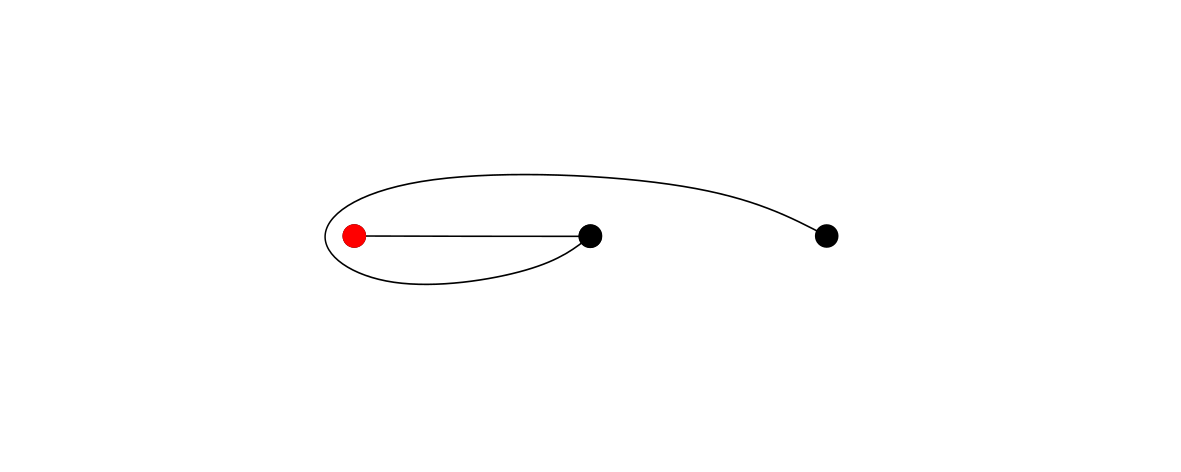}\par 
			\includegraphics[width=\linewidth]{ex7.png}\par 
		\end{multicols}
		\vspace{-1 cm} 
		\begin{multicols}{3}
			\hspace{.1 cm}\par 
			\makebox[\textwidth][l]{\textbf{step 3:} \hspace{0 cm} \text{$ g_1' =  s_2^{-1}.s_2^{-1}$}}\par 
			\makebox[\textwidth][l]{\textbf{step 3:} \hspace{0 cm} \text{$ g_2' = id$}}\par 
		\end{multicols}
		\vspace{-1 cm} 
		\begin{multicols}{3}
			\hspace{.1 cm} \par   
			\includegraphics[width=\linewidth]{ex6.png}\par
			\includegraphics[width=\linewidth]{ex6.png}\par 
		\end{multicols}
		\vspace{-1 cm} 
		\begin{multicols}{3}
			\hspace{.1 cm}\par 
			\makebox[\textwidth][l]{\textbf{step 4:} \hspace{.0 cm} \text{$ g_1 =  s_2^{-1}.s_2^{-1}$}}\par 
			\makebox[\textwidth][l]{\textbf{step 4:} \hspace{0 cm} \text{$ g_2 = s_2.s_1.s_0.s_0.s_1.s_2$}}\par 
		\end{multicols}
		\vspace{-1 cm} 
		\begin{multicols}{3}
			\hspace{.1 cm} \par   
			\includegraphics[width=\linewidth]{ex6.png}\par
			\includegraphics[width=\linewidth]{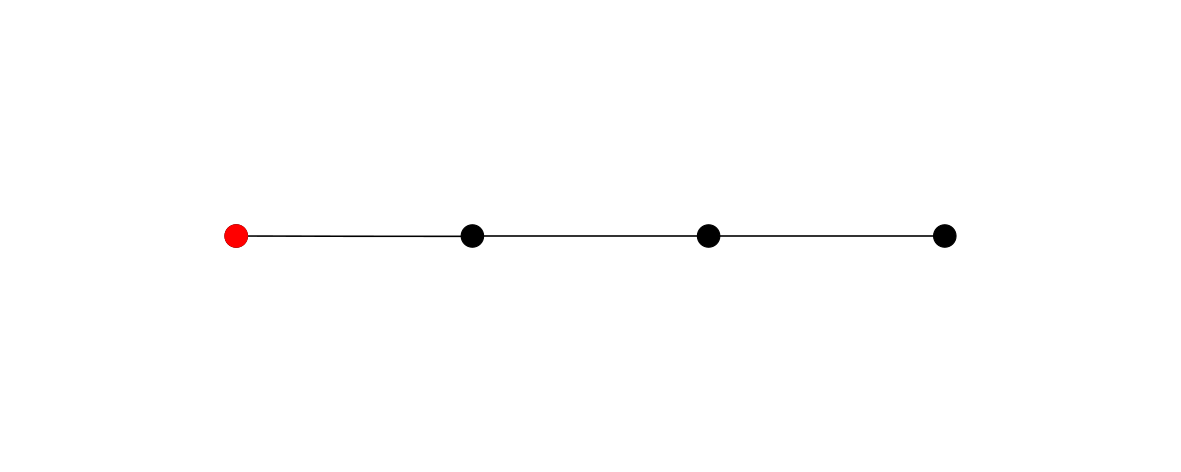} \par
		\end{multicols}	
			\caption{The algorithm described in Section \ref{sec:algorithm} is applied in a simple case. We obtain $g = g_2\circ g_1\circ g_0 =  s_2.s_1.s_0.s_0.s_1.s_2.s_{2}^{-1}.s_{2}^{-1}.s_1$.}
		\label{fig:drawing}	
	\end{figure}

\begin{figure}
		\centering
		\vspace{-.25 cm}
		\hspace{2 cm }\includegraphics[width=0.6\textwidth]{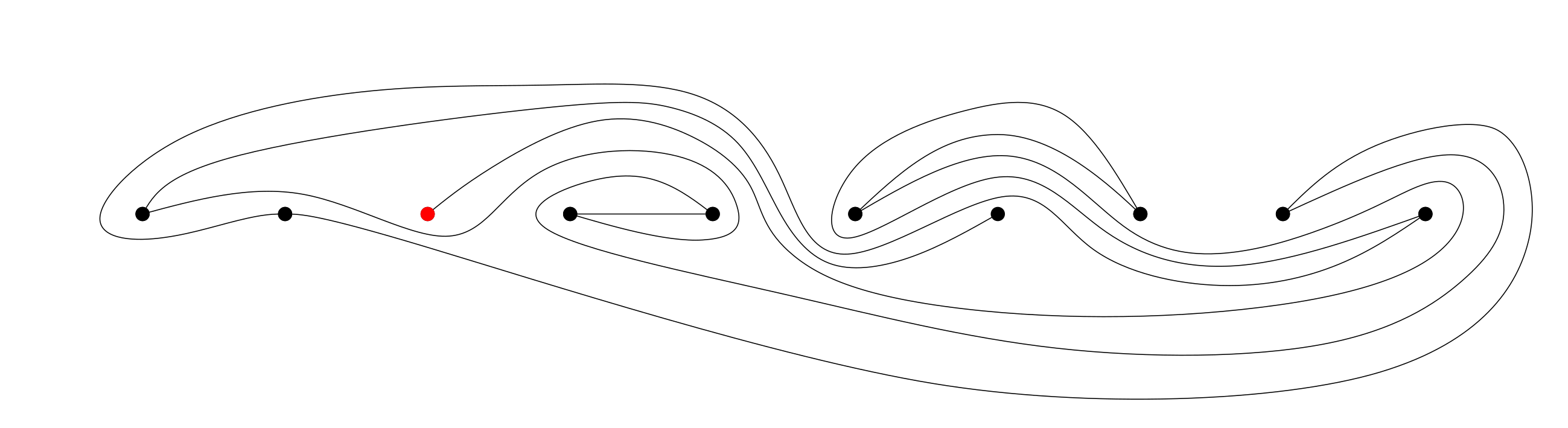} 
		\makebox[\textwidth][l]{\textbf{apply: $g_0 = $} \hspace{3.1 cm} $s_1.s_2.s_0.s_1.s_3.s_4.s_5^{-1}.s_6^{-1}.s_7^{-1}.s_8^{-1}.s_8^{-1}.s_7.s_6.s_5^{-1}$}
		\\ \hspace{2 cm} \includegraphics[width=0.6\textwidth]{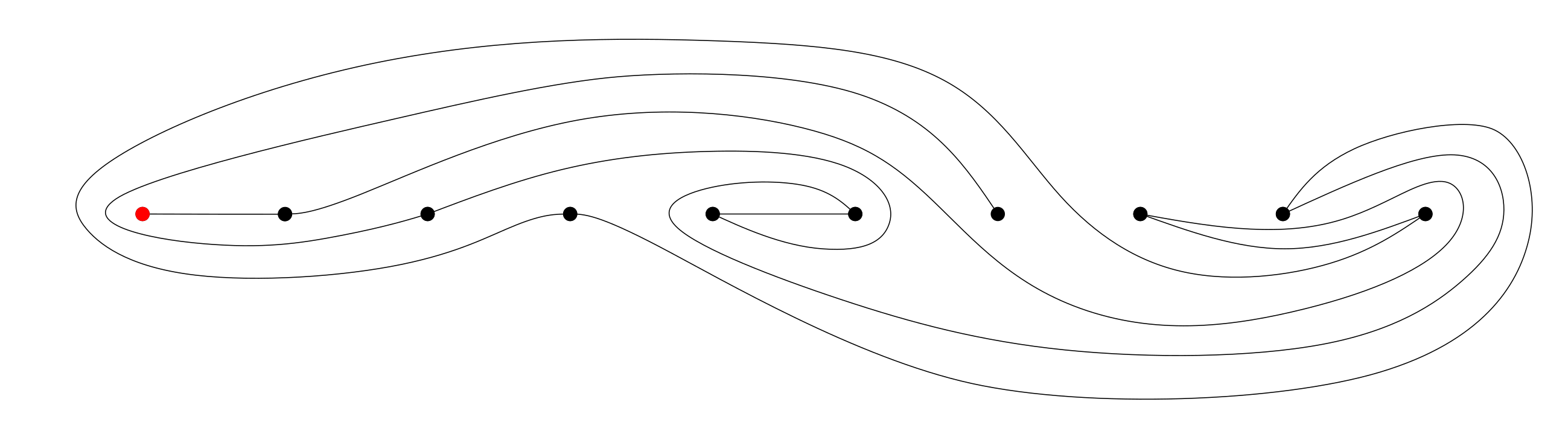}
		\makebox[\textwidth][l]{\textbf{apply: $g_1 = $} \hspace{4.5 cm} $s_2.s_3.s_4.s_5.s_6^{-1}.s_7^{-1}.s_8^{-1}.s_8^{-1}.s_7$}
		\\ \hspace{2 cm} \includegraphics[width=0.6\textwidth]{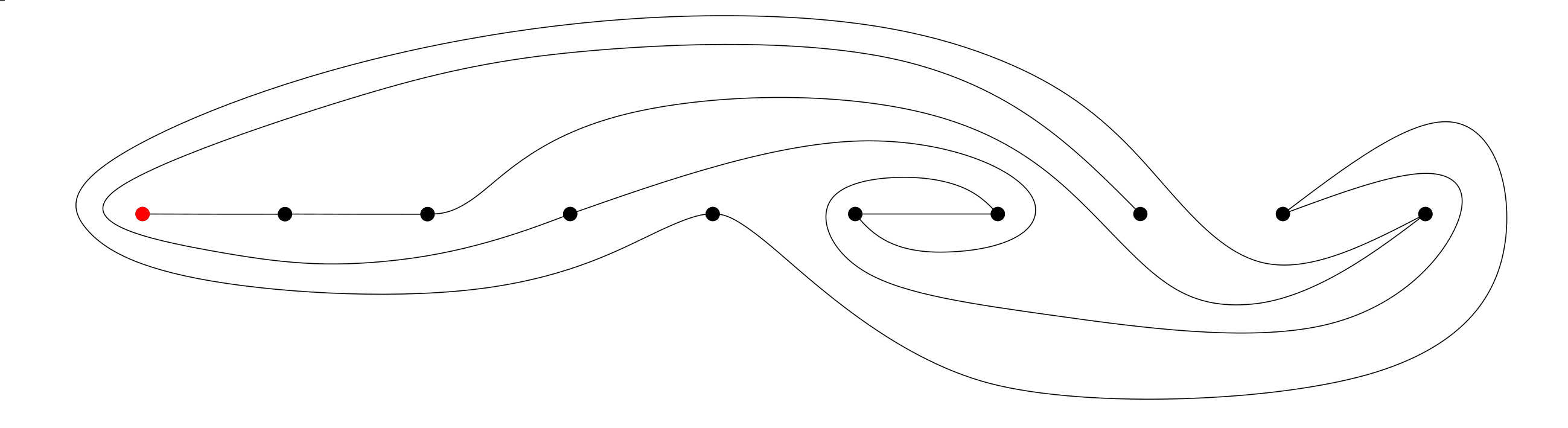}
		\makebox[\textwidth][l]{\textbf{apply: $g_2 = $} \hspace{5.3 cm} $s_3.s_4.s_5.s_6.s_7^{-1}.s_8^{-1}$}   
		\\ \hspace{2 cm} \includegraphics[width=0.6\textwidth]{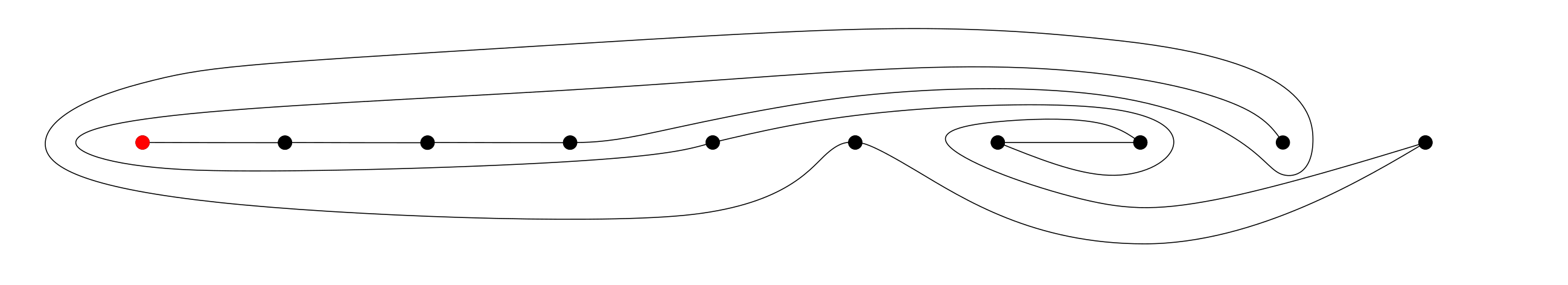}
		\makebox[\textwidth][l]{\textbf{apply: $g_3 = $} \hspace{2 cm} $s_4.s_5.s_6.s_7^{-1}.s_7^{-1}.s_6^{-1}.s_5^{-1}.s_4^{-1}.s_3^{-1}.s_2^{-1}.s_1^{-1}.s_0^{-1}.s_0^{-1}.s_1^{-1}.s_2^{-1}.s_3^{-1}.s_4^{-1}$} 
		\\ \hspace{2 cm}  \includegraphics[width=0.6\textwidth]{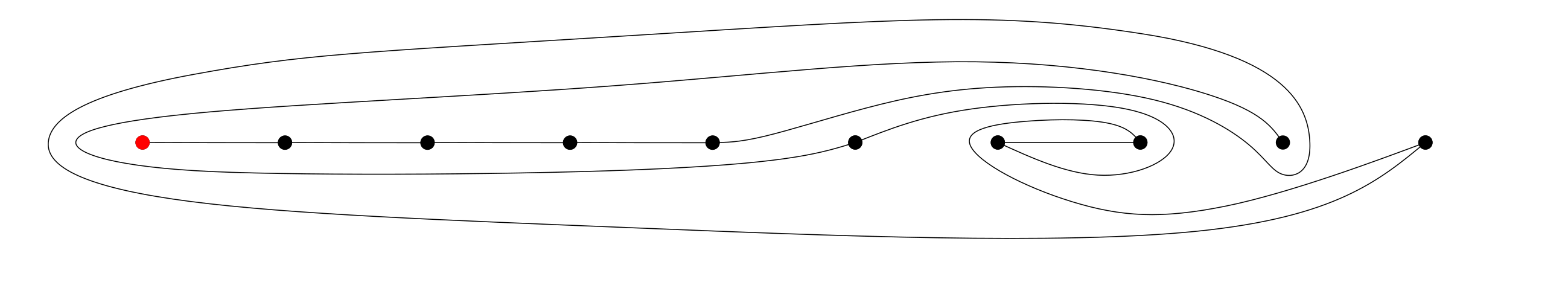}
		\makebox[\textwidth][l]{\textbf{apply: $g_4 = $} \hspace{.5  cm} $s_5.s_6.s_7.s_8^{-1}.s_8^{-1}.s_7^{-1}.s_6^{-1}.s_5^{-1}.s_4^{-1}.s_3^{-1}.s_2^{-1}.s_1^{-1}.s_0^{-1}.s_0^{-1}.s_1^{-1}.s_2^{-1}.s_3^{-1}.s_4^{-1}.s_5^{-1}.s_6^{-1}.s_7^{-1}.s_8^{-1}$}
		\\ \hspace{2 cm} \includegraphics[width=0.6\textwidth]{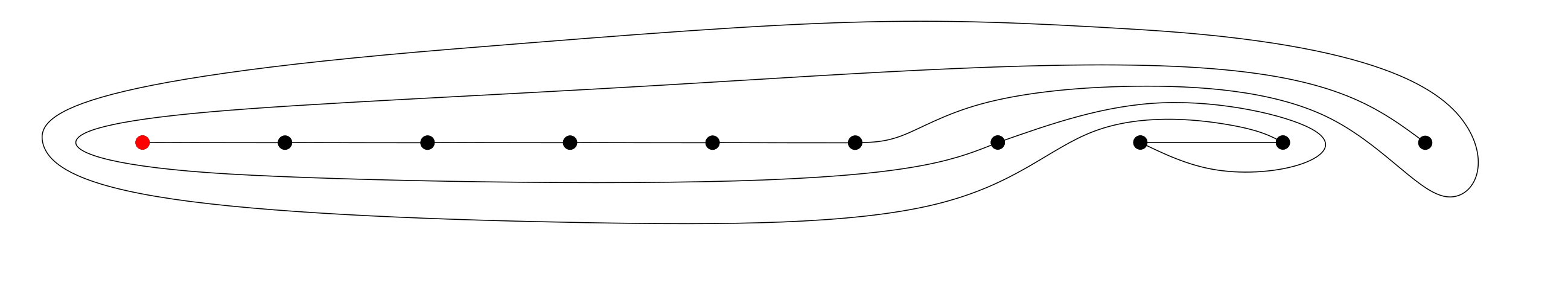} 
		\makebox[\textwidth][l]{\textbf{apply: $g_5 = $} \hspace{1 cm} $s_6.s_7.s_8^{-1}.s_8^{-1}.s_7^{-1}.s_6^{-1}.s_5^{-1}.s_4^{-1}.s_3^{-1}.s_2^{-1}.s_1^{-1}.s_0^{-1}.s_0^{-1}.s_1^{-1}.s_2^{-1}.s_3^{-1}.s_4^{-1}.s_5^{-1}.s_6^{-1}.s_7$}    
		\\ \hspace{2 cm} \includegraphics[width=0.6\textwidth]{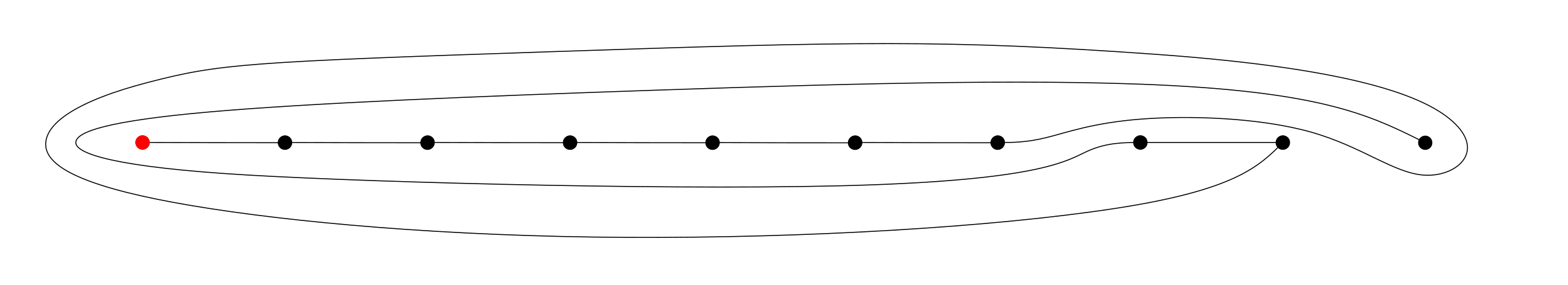}
		\makebox[\textwidth][l]{\textbf{apply: $g_6 = $} \hspace{1.1 cm}  $s_7.s_8^{-1}.s_8^{-1}.s_7^{-1}.s_6^{-1}.s_5^{-1}.s_4^{-1}.s_3^{-1}.s_2^{-1}.s_1^{-1}.s_0^{-1}.s_0^{-1}.s_1^{-1}.s_2^{-1}.s_3^{-1}.s_4^{-1}.s_5^{-1}.s_6^{-1}.s_7^{-1}$} 
		\\ \hspace{2 cm} \includegraphics[width=0.6\textwidth]{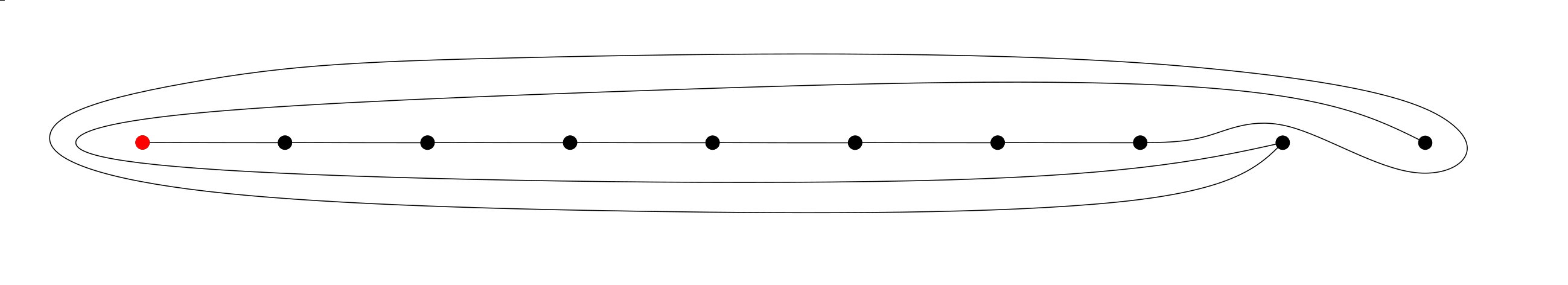}
		\makebox[\textwidth][l]{\textbf{apply: $g_7 = $} \hspace{1.2cm} $s_8^{-1}.s_8^{-1}.s_7^{-1}.s_6^{-1}.s_5^{-1}.s_4^{-1}.s_3^{-1}.s_2^{-1}.s_1^{-1}.s_0^{-1}.s_0^{-1}.s_1^{-1}.s_2^{-1}.s_3^{-1}.s_4^{-1}.s_5^{-1}.s_6^{-1}.s_7^{-1}$} 
		\\ \hspace{2 cm} \includegraphics[width=0.6\textwidth]{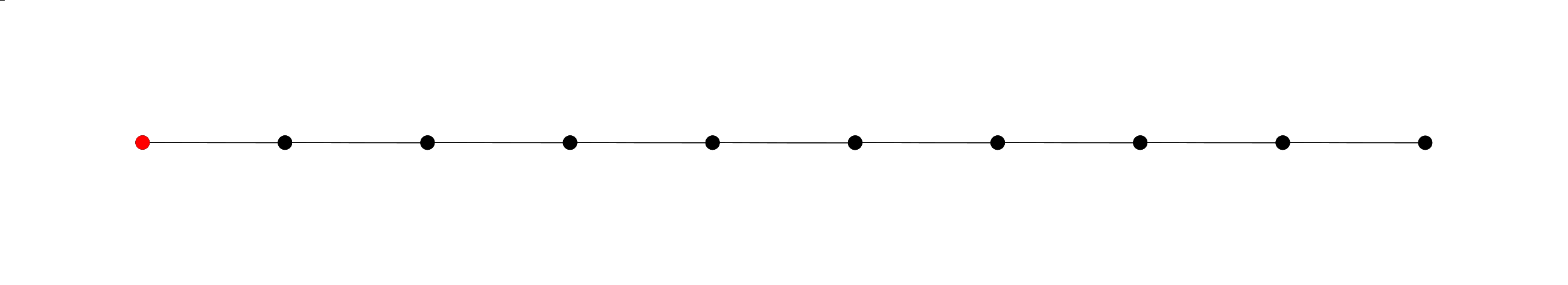}
		\vspace{-.5 cm} 
		\caption{The top panel depicts representatives of the relative homotopy classes of the paths $f'\circ p_i$, as certified in Section \ref{sec:arc-certification}. The rest of the figure applies the algorithm described in Section \ref{sec:algorithm} to the mapping class induced by $f$ on the sphere punctured along the $10$-periodic orbit obtained in Section \ref{sec:ApplyShadowing}.}
		\label{fig:actual}
	\end{figure}

	For all $0 \leq i < n-1$, we will find $g_i \in \operatorname{Homeo}^+(\mathring{D}_n)$ such that $g_i(p_k) \cong p_k$ for each $0\leq k < i$ and $g_i\circ g_{i-1}\circ \dots \circ g_0\circ f'\circ p_i \cong p_i$. Then 
	$g \coloneqq g_{n-2}\circ g_{n-3}\circ \dots \circ g_0$ satisfies $g\circ f'\circ p_i\cong p_i$ for each $0\leq i < n-1$, as desired.
	We will determine $g_i$ inductively in the following way.
	\vspace{.25 cm} \\
	\noindent	{\bf Base case:}  Choose a finite sequence of half-twists $t_1, \dots, t_\ell\in \set{s_k^\pm \mid 0\leq k< n-1}$ such that 
	\begin{align*}
			g_0 \coloneqq t_\ell.t_{\ell-1}.\dots . t_1\in  \operatorname{Homeo}^+(\mathring D_{n})
		\end{align*} satisfies $g_0\circ f'\circ p_0 \cong p_0$.
	The half-twists $t_k$ can be readily obtained from the arc-data of a representative of the relative homotopy class of $f'\circ p_0$. The process is straightforward and we omit details; one can follow along in Figures \ref{fig:drawing} and \ref{fig:actual} to get a sense of it. The idea is to ``drag" the endpoint of the representative backwards along the path towards its starting point using half-twists.
	{\color{orange}\leavevmode}
	\vspace{.25 cm} \\
	\noindent	{\bf Inductive case:} Let $i > 0$, and assume that $g_k$ for $0\leq k  < i$ have been constructed.
	\vspace{.25 cm}\\
	\indent {\bf Step 1.}  Define $P_{i-1} \coloneqq p_0.p_1. \dots .p_{i-1}$. Define a homotopy $C: [0,1]\times \mathring{D} \to \mathring{D}$ such that:
	\begin{itemize}
		\item $C_0$ is the identity;
		\item $C_1(P_{i-1})$ is a point;
		\item $C_1$ restricts to a homeomorphism from $\mathring D\setminus P_{i-1}$ to $\mathring D\setminus\set{C_1(P_{i-1})}$;
		\item for all $t\in [0,1)$, $C_t$ is a homeomorphism satisfying $C_t(P_{i-1})\subset P_{i-1}$;
		\item for all $t\in [0,1]$, $C_t$ fixes $z_k$ whenever $k> i$.
	\end{itemize}
	Let $c\coloneqq C_1$. Considering $c(P_{i-1})$ to be a marked point along with the $z_k$'s for $i < k \leq n-1$, $c$ defines a continuous map $c: \mathring{D}_n \to \mathring{D}_{n - i}$.
		We take generators $\set{\mathfrak{s}_k \mid i\leq k< n-1}$ for $\MOD^+(\mathring{D}_{n-i})$, where $\mathfrak{s}_k$ is a counter-clockwise half-twist about $c(p_k)$. Notice that, for $i < k < n-1$, one has 
	\begin{equation}
		\mathfrak{s}_k\circ c \cong c \circ s_k. \label{lift0}
	\end{equation} Moreover, the full twist $\tau_{c(p_i)}=\mathfrak{s}_i^2$ admits the following lift by $c$:
	\begin{align}
		\tau_{c(p_i)}\circ c
		&\cong c\circ (s_i.s_{i-1}.\dots.s_1.s_0.s_0.s_1.\dots.s_{i-1}.s_i).
		\label{lift}
	\end{align}
	
		{\bf Step 2.} Consider the path $\alpha_i\coloneqq g_{i-1}\circ g_{i-2}\circ \dots \circ g_0\circ f'\circ p_i:[0,1]\to \mathring{D}_n$ and 
	$c(\alpha_i)\subset \mathring{D}_{n-i}$. Choose $\tilde g_i\in \operatorname{Homeo}^+(\mathring D_{n-i})$ such that $\tilde g_i(c(\alpha_i)) \cong c(p_i)$ in $\mathring D_{n-i}$. Notice that $c(\alpha_i)$ and $c(p_i)$ both have starting point $c(P_{i-1})$. Therefore, $ [\tilde g_i]\in \MOD^+_{0}(\mathring{D}_{n-i})$, and by Lemma \ref{lem:generators}, $ [\tilde g_i ]\in \angle{[\tau_{c(p_i)}],[ \mathfrak{s}_{i+1}], [\mathfrak{s}_{i+2}], \dots,[ \mathfrak{s}_{n-2}]}$.
	Just like in the base case, a presentation for $[\tilde g_i]$ in terms of those generators can be obtained from the arc-data of a representative of the relative homotopy class of $c\circ\alpha_i$.
	
	{\bf Step 3.}
	By applying \eqref{lift0} and  \eqref{lift}, we can
	lift $\tilde g_i$ to $g_i'\in \operatorname{Homeo}^+(\mathring{D}_n)$ satisfying $\tilde g_i\circ c\cong c\circ g_i'$.
	The homeomorphism $ g_i '$ fixes $p_k$ whenever $0\leq k < i$.
	Moreover, $g_i'$ maps $\alpha_i$ to $p_i$ up to twisting around $P_{i-1}$; we make this precise in the next step.
	
		{\bf Step 4.}

	Let $P_i\coloneqq P_{i-1}.p_i$, and let $\tau_{P_{i-1}}$ and $\tau_{P_i}$ denote the counter-clockwise Dehn twists about $P_{i-1}$ and $P_i$, respectively.
	Choose representatives so that $P_{i-1}.g_i'(\alpha_i)$ is embedded. Then $g_i'(\alpha_i)$ and $p_i$ have endpoints $z_i$ and $z_{i+1}$, and
	\begin{align*}
		c(g_i'(\alpha_i))\cong\tilde g_i(c(\alpha_i))\cong c(p_i).
	\end{align*}
	We use the following elementary fact.
	\begin{claim}\label{claim:gen2}
		Suppose that $\beta$ is an embedded arc from $z_i$ to $z_{i+1}$ satisfying $\beta\cap P_{i-1}=\set{z_i}$ and $c(\beta)\cong c(p_i)$ relative to the marked points of the collapsed disk. Then there exists $k\in\Z$ such that
		\begin{align*}
			p_i\cong\tau_{P_{i-1}}^k(\beta).
		\end{align*}
	\end{claim}
	\begin{proof}[Proof of Claim \ref{claim:gen2}]
		Remove a small regular neighborhood of $P_{i-1}$. The arcs $\beta$ and $p_i$ become proper arcs beginning on the same boundary component. Since their images under $c$ are homotopic, their only possible difference is an integer winding about this boundary component. This winding is generated by $\tau_{P_{i-1}}$, proving the claim.
	\end{proof}
	Applying Claim \ref{claim:gen2} with $\beta=g_i'(\alpha_i)$, we find that
	\begin{align*}
		p_i\cong\tau_{P_{i-1}}^k(g_i'(\alpha_i))
	\end{align*}
	for some $k\in\Z$.
	One could let $g_i \coloneqq \tau_{P_{i-1}}^k\circ  g_i'$, but we instead set 
	$g_i \coloneqq  \tau_{P_{i}}^{-k}\circ \tau_{P_{i-1}}^k\circ g_i'$. Notice that $\tau_{P_{i}}$ and  $\tau_{P_{i-1}}$ commute; this can be seen by realizing $\tau_{P_{i}}$ and $\tau_{P_{i-1}}$ as twists around disjoint closed loops.
	So 
	\begin{align*}
		\tau_{P_i}^{-k}\tau_{P_{i-1}}^k
		&=(\tau_{P_i}^{-1}\tau_{P_{i-1}})^k
			=((s_i^{-1}.\dots.s_0^{-1}).(s_0^{-1}.\dots.s_i^{-1}))^k,
	\end{align*}
	and we have
	\begin{align*}
		g_i =((s_i^{-1}.\dots.s_0^{-1}).(s_0^{-1}.\dots.s_i^{-1}))^k\circ g_i',
	\end{align*} concluding the inductive step.

	\subsubsection{Application to our case and use of Flipper} \label{section:application}
	In Figure \ref{fig:actual}, we apply the algorithm described in Section \ref{sec:algorithm} to compute $g$. After applying $g_7\circ\dots\circ g_0$, we have
	\begin{align*}
		g_7\circ\dots\circ g_0\circ f'\circ p_8\cong p_8,
	\end{align*}
	so we set $g_8=\operatorname{id}$. We find:
	\begin{align*}
		g &= g_8\circ g_{7}\circ g_{6}\circ \dots \circ g_0 = g_{7}\circ g_{6}\circ \dots \circ g_0 \\ 
		&= (s_8^{-1}.s_8^{-1}.s_7^{-1}.s_6^{-1}.s_5^{-1}.s_4^{-1}.s_3^{-1}.s_2^{-1}.s_1^{-1}.s_0^{-1}.s_0^{-1}.s_1^{-1}.s_2^{-1}.s_3^{-1}.s_4^{-1}.s_5^{-1}.s_6^{-1}.s_7^{-1})\\
		&\hspace{.75 cm} \circ
		(s_7.s_8^{-1}.s_8^{-1}.s_7^{-1}.s_6^{-1}.s_5^{-1}.s_4^{-1}.s_3^{-1}.s_2^{-1}.s_1^{-1}.s_0^{-1}.s_0^{-1}.s_1^{-1}.s_2^{-1}.s_3^{-1}.s_4^{-1}.s_5^{-1}.s_6^{-1}.s_7^{-1}) \\
		&\hspace{.75 cm} \circ (s_6.s_7.s_8^{-1}.s_8^{-1}.s_7^{-1}.s_6^{-1}.s_5^{-1}.s_4^{-1}.s_3^{-1}.s_2^{-1}.s_1^{-1}.s_0^{-1}.s_0^{-1}.s_1^{-1}.s_2^{-1}.s_3^{-1}.s_4^{-1}.s_5^{-1}.s_6^{-1}.s_7) \\
		&\hspace{.75 cm} \circ  (s_5.s_6.s_7.s_8^{-1}.s_8^{-1}.s_7^{-1}.s_6^{-1}.s_5^{-1}.s_4^{-1}.s_3^{-1}.s_2^{-1}.s_1^{-1}.s_0^{-1}.s_0^{-1}.s_1^{-1}.s_2^{-1}.s_3^{-1}.s_4^{-1}.s_5^{-1}.s_6^{-1}.s_7^{-1}.s_8^{-1})\\
		& \hspace{.75 cm}  \circ (s_4.s_5.s_6.s_7^{-1}.s_7^{-1}.s_6^{-1}.s_5^{-1}.s_4^{-1}.s_3^{-1}.s_2^{-1}.s_1^{-1}.s_0^{-1}.s_0^{-1}.s_1^{-1}.s_2^{-1}.s_3^{-1}.s_4^{-1})  \circ ( s_3.s_4.s_5.s_6.s_7^{-1}.s_8^{-1}) \\
		&\hspace{.75 cm} 
		\circ
		(s_2.s_3.s_4.s_5.s_6^{-1}.s_7^{-1}.s_8^{-1}.s_8^{-1}.s_7)\circ
		(s_1.s_2.s_0.s_1.s_3.s_4.s_5^{-1}.s_6^{-1}.s_7^{-1}.s_8^{-1}.s_8^{-1}.s_7.s_6.s_5^{-1}).
	\end{align*}
	As we saw in \eqref{fSquared}, $[f^2]$ is readily obtained from $g$. For completeness, we include its mapping class below:
\small
	\begin{equation}\label{eq:f2-flipper-word}
			\begin{aligned}
			f^2 &\cong
			\big[(\hat s_5.\hat s_6^{-1}.\hat s_7^{-1}.\hat s_8.\hat s_8.\hat s_7.\hat s_6.\hat s_5.\hat s_4^{-1}.\hat s_3^{-1}.\hat s_1^{-1}.\hat s_0^{-1}.\hat s_2^{-1}.\hat s_1^{-1})
			\circ
			(\hat s_7^{-1}.\hat s_8.\hat s_8.\hat s_7.\hat s_6.\hat s_5^{-1}.\hat s_4^{-1}.\hat s_3^{-1}.\hat s_2^{-1})	\\
			&\hspace{.75 cm} 
			\circ (\hat s_8.\hat s_7.\hat s_6^{-1}.\hat s_5^{-1}.\hat s_4^{-1}.\hat s_3^{-1})
			\circ	(\hat s_4.\hat s_3.\hat s_2.\hat s_1.\hat s_0.\hat s_0.\hat s_1.\hat s_2.\hat s_3.\hat s_4.\hat s_5.\hat s_6.\hat s_7.\hat s_7.\hat s_6^{-1}.\hat s_5^{-1}.\hat s_4^{-1})\\
			&\hspace{.75 cm} \circ
			(\hat s_8.\hat s_7.\hat s_6.\hat s_5.\hat s_4.\hat s_3.\hat s_2.\hat s_1.\hat s_0.\hat s_0.\hat s_1.\hat s_2.\hat s_3.\hat s_4.\hat s_5.\hat s_6.\hat s_7.\hat s_8.\hat s_8.\hat s_7^{-1}.\hat s_6^{-1}.\hat s_5^{-1})\\
			&\hspace{.75 cm} \circ
			(	\hat s_7^{-1}.\hat s_6.\hat s_5.\hat s_4.\hat s_3.\hat s_2.\hat s_1.\hat s_0.\hat s_0.\hat s_1.\hat s_2.\hat s_3.\hat s_4.\hat s_5.\hat s_6.\hat s_7.\hat s_8.\hat s_8.\hat s_7^{-1}.\hat s_6^{-1})\\
			&\hspace{.75 cm} \circ
			(\hat s_7.\hat s_6.\hat s_5.\hat s_4.\hat s_3.\hat s_2.\hat s_1.\hat s_0.\hat s_0.\hat s_1.\hat s_2.\hat s_3.\hat s_4.\hat s_5.\hat s_6.\hat s_7.\hat s_8.\hat s_8.\hat s_7^{-1})	\\
			&\hspace{.75 cm} 
			\circ (\hat s_7.\hat s_6.\hat s_5.\hat s_4.\hat s_3.\hat s_2.\hat s_1.\hat s_0.\hat s_0.\hat s_1.\hat s_2.\hat s_3.\hat s_4.\hat s_5.\hat s_6.\hat s_7.\hat s_8.\hat s_8)\big] \\
			&\hspace{.25 cm} 
			\circ \big[(\hat s_5^{-1}.\hat s_6.\hat s_7.\hat s_8^{-1}.\hat s_8^{-1}.\hat s_7^{-1}.\hat s_6^{-1}.\hat s_5^{-1}.\hat s_4.\hat s_3.\hat s_1.\hat s_0.\hat s_2.\hat s_1)
			\circ (\hat s_7.\hat s_8^{-1}.\hat s_8^{-1}.\hat s_7^{-1}.\hat s_6^{-1}.\hat s_5.\hat s_4.\hat s_3.\hat s_2) \\
			&\hspace{.75 cm} 
			\circ (\hat s_8^{-1}.\hat s_7^{-1}.\hat s_6.\hat s_5.\hat s_4.\hat s_3) 
			\circ (\hat s_4^{-1}.\hat s_3^{-1}.\hat s_2^{-1}.\hat s_1^{-1}.\hat s_0^{-1}.\hat s_0^{-1}.\hat s_1^{-1}.\hat s_2^{-1}.\hat s_3^{-1}.\hat s_4^{-1}.\hat s_5^{-1}.\hat s_6^{-1}.\hat s_7^{-1}.\hat s_7^{-1}.\hat s_6.\hat s_5.\hat s_4) \\
			&\hspace{.75 cm} 
			\circ (\hat s_8^{-1}.\hat s_7^{-1}.\hat s_6^{-1}.\hat s_5^{-1}.\hat s_4^{-1}.\hat s_3^{-1}.\hat s_2^{-1}.\hat s_1^{-1}.\hat s_0^{-1}.\hat s_0^{-1}.\hat s_1^{-1}.\hat s_2^{-1}.\hat s_3^{-1}.\hat s_4^{-1}.\hat s_5^{-1}.\hat s_6^{-1}.\hat s_7^{-1}.\hat s_8^{-1}.\hat s_8^{-1}.\hat s_7.\hat s_6.\hat s_5) \\
			&\hspace{.75 cm} 
			\circ (\hat s_7.\hat s_6^{-1}.\hat s_5^{-1}.\hat s_4^{-1}.\hat s_3^{-1}.\hat s_2^{-1}.\hat s_1^{-1}.\hat s_0^{-1}.\hat s_0^{-1}.\hat s_1^{-1}.\hat s_2^{-1}.\hat s_3^{-1}.\hat s_4^{-1}.\hat s_5^{-1}.\hat s_6^{-1}.\hat s_7^{-1}.\hat s_8^{-1}.\hat s_8^{-1}.\hat s_7.\hat s_6) \\
			&\hspace{.75 cm} 
				\circ(\hat s_7^{-1}.\hat s_6^{-1}.\hat s_5^{-1}.\hat s_4^{-1}.\hat s_3^{-1}.\hat s_2^{-1}.\hat s_1^{-1}.\hat s_0^{-1}.\hat s_0^{-1}.\hat s_1^{-1}.\hat s_2^{-1}.\hat s_3^{-1}.\hat s_4^{-1}.\hat s_5^{-1}.\hat s_6^{-1}.\hat s_7^{-1}.\hat s_8^{-1}.\hat s_8^{-1}.\hat s_7) \\
				&\hspace{.75 cm}
				\circ(\hat s_7^{-1}.\hat s_6^{-1}.\hat s_5^{-1}.\hat s_4^{-1}.\hat s_3^{-1}.\hat s_2^{-1}.\hat s_1^{-1}.\hat s_0^{-1}.\hat s_0^{-1}.\hat s_1^{-1}.\hat s_2^{-1}.\hat s_3^{-1}.\hat s_4^{-1}.\hat s_5^{-1}.\hat s_6^{-1}.\hat s_7^{-1}.\hat s_8^{-1}.\hat s_8^{-1})\big].
		\end{aligned}
	\end{equation}

	\normalsize
	
			In \href{https://github.com/ethanhcoo/K3entropy.git}{K3entropy}, \texttt{FlipperDilatation.ipynb}\footnotemark{} proves the following Proposition. Recall the polynomial $P$ defined in \eqref{definedP}.
 \footnotetext{Flipper uses right Dehn twists. The distinction between left
	and right is not important here, since left and right twists are conjugate
	by an orientation-reversing homeomorphism, and conjugation preserves
	topological entropy.}
	\begin{prop}[Flipper certificate]\label{prop:FlipperCertificate}
		The mapping class $[f^2]$ is pseudo-Anosov, and its stretching factor is the unique real root of $P$ in $(8.19981448, 8.19981449)$, which we denote by $\lambda$. 
	\end{prop}
	
	\begin{proof}
		Flipper's Nielsen--Thurston classification routine verifies that the mapping class represented by the displayed word is pseudo-Anosov. It then constructs the $42\times42$ hitting matrix $M$ for its invariant train track. The full matrix is printed by
		\texttt{FlipperDilatation.ipynb}, and its characteristic polynomial factors as
		\[\chi_M(t) =P(t)(t-1)^4(t+1)^4(t^2+t+1) (t^{18}+t^{17}+\cdots+t+1).\]
		 Then SageMath verifies that $P$ is irreducible, and exact real-root isolation isolates all of its positive real roots. The largest is the unique root $\lambda$ in $(8.19981448,8.19981449)$, and $P$ changes sign across the endpoints. Since $M$ is nonnegative, its spectral radius is a positive real eigenvalue, and all the other factors in $\chi_M$ are cyclotomic. Therefore, $\lambda$ is the spectral radius of $M$, and its minimal polynomial is $P$. This proves the claim.
	\end{proof}
	
	\noindent It follows that
	\[h_{top}(f, X(\R))\geq \frac{1}{2}\ln\lambda. \]
	Combining this with \eqref{ComplexEntropy} gives
	\begin{align*}
		\frac{h_{top}(f, X(\R))}{h_{top}(f, X(\C))}
		\geq \frac{\ln\lambda}{2\ln\mu} \approx 0.58
		> \frac 1 2.
	\end{align*}

	\begin{appendices}

		\section{Proof of the shadowing lemma}\label{app:closing}
		\begin{proof}[Proof of Lemma \ref{lem:closing}]
			To start, we set up a fixed point problem. Let $\epsilon \coloneqq 12\kappa\delta$ so that $B_{\epsilon}(0)\subset U_i$.
			Let
			\begin{align*}
				K\coloneqq\max_{0\leq i<n}\norm{D^2h_i^c}_{C^0(B_{\epsilon}(0)),\infty}.
			\end{align*}
			Pick a smooth bump function $g:\R^2\to \R$ satisfying $0\leq g\leq1$, $\operatorname{supp} g\subset B_{\epsilon}(0)$, and $\restr{g}{B_{\frac{\epsilon}{2}} (0)}\equiv 1$. Furthermore, choose $g$ such that $\norm{Dg}_{C^0(\R^2)} < \frac{3}{\epsilon}$. For each $0\leq i<n$, let
			\begin{align*}
				L_i\coloneqq(Dh^c_i)_0.
			\end{align*}
			On $U_i$, set
			\begin{align*}
				H_i\coloneqq g h^c_i + (1 - g)L_i.
			\end{align*}
			Since $\operatorname{supp} g\subset B_{\epsilon}(0)\subset U_i$, this extends to a $C^2$ map $H_i:\R^2\to\R^2$ by setting $H_i=L_i$ outside $U_i$. Finally, let $M_i\coloneqq H_i-L_i:\R^2\to\R^2$.
			\noindent Furthermore, define $H, L, M:\R^{2n}\to \R^{2n}$ by 
			\begin{align*}
				H(w_0, \dots, w_{n-1})&\coloneqq (H_{n-1}(w_{n-1}), H_0(w_0),\dots, H_{n-2}(w_{n-2})),\\
				L(w_0, \dots, w_{n-1})&\coloneqq (L_{n-1}(w_{n-1}), L_0(w_0),\dots, L_{n-2}(w_{n-2})), \\
				M(w_0, \dots, w_{n-1})&\coloneqq (M_{n-1}(w_{n-1}), M_0(w_0),\dots, M_{n-2}(w_{n-2})).
			\end{align*}
			We hope to find a fixed point $z = (z_0, \dots, z_{n-1})$ for $H$ in $B_{\frac{\epsilon}{2}}(\vec 0)$. 
			Here we are using the supremum norm $\norm{\cdot}_0$ on $(\R^2)^n$ (see \ref{sec:notations}). 
			Notice that we require $z_i\in B_{\frac{\epsilon}{2}}(0)$ for every $0\leq i<n$ in order to lift $z$ to an $h$-periodic point. Indeed, since $\restr{H_i}{B_{\frac {\epsilon}{2}}(0)} \equiv h_i^c$, one has $H(z) = (h_{n-1}^c(z_{n-1}), h_0^c(z_0), \dots, h_{n-2}^c(z_{n-2}))$. Thus  $(\tilde x_i)\coloneqq (\phi_i(z_i))$ is a fixed point for $(w_0, \dots, w_{n-1})\mapsto (h(w_{n-1}), h(w_0),\dots, h(w_{n-2}))$ as a function $M^n\to M^n$.
			Note that 
			\begin{align*}
				z = H(z) = L(z) + M(z)
			\end{align*} exactly when
			\begin{align*}
				(L - id)z =  -M(z)
			\end{align*} or in other words
			\begin{align*}
				z =  -(L - id)^{-1}M(z).
			\end{align*} So, we hope to find a fixed point for $\mathcal{F}\coloneqq -(L - id)^{-1}M$ in $B_{\frac{ \epsilon}{2}}(\vec 0)$. 
			
			\begin{claim}\label{claim:first}
				$\mathcal{F}$ is a $\frac{13}{16}$-contraction mapping for $\norm{\cdot}_0$ whenever 
				\[\epsilon < \frac{1}{16\kappa K}.\]
			\end{claim}
			\noindent We must prove $\sup_{w\in \R^{2n}}\norm{D((L - id)^{-1}M)_w}_{op,0} \leq \frac{13}{16}$. We treat $\sup_{w\in \R^{2n}}\norm{D((L - id)^{-1})_w}_{op,0}$ and  $\sup_{w\in \R^{2n}}\norm{DM}_{op,0}$ separately. Now $L$ is linear, so $D((L - id)^{-1})_w = (L - id)^{-1}$ for all $w$. Thus $\sup_{w\in \R^{2n}}\norm{D((L - id)^{-1})_w}_{op,0} = \kappa$. 
			
			\begin{claim}\label{claim:second}
					$\sup_{w\in \R^{2n}}\norm{(DM)_w}_{op,0}
					\leq \frac{3\delta}{\epsilon}+(6+2\sqrt{2})\epsilon K
					<\frac{3\delta}{\epsilon}+9\epsilon K.$
			\end{claim}
			\noindent Assuming \ref{claim:second}, we have
			\begin{align*}
				\sup_{w\in \R^{2n}}\norm{D((L - id)^{-1}M)_w}_{op,0}
				&<\kappa\left(9\cdot\frac{1}{16\kappa}+\frac{3}{12\kappa}\right) =\frac{13}{16},
			\end{align*} proving \ref{claim:first}. 
			\noindent To begin the proof of \ref{claim:second}, define $N_i^c:B_{\epsilon}(0)\to \R^2$ by
			\begin{align*}
				N_i^c &\coloneqq  h_i^c - L_i.
			\end{align*} The function $N_i^c$ is essentially the non-linear part of $h^c_i$ at $0$, with the caveat that $N_i^c$ might not fix the origin. On $B_{\epsilon}(0)$, one has $M_i = H_i - L_i =  g(h_i^c - L_i) = g N^c_i$, while $M_i$ vanishes outside $\operatorname{supp} g$. Therefore,
			\begin{align*}
				\sup_{w\in \R^2}\norm{(DM_i)_w}_{op}& \leq \norm{DM_i}_{C^0(\R^2)} \\ & \leq \norm{Dg}_{C^0(\R^2)} \norm{N_i^c}_{C^0(B_{\epsilon}(0))}  + \norm{g}_{C^0(\R^2)}\norm{DN_i^c}_{C^0(B_{\epsilon}(0))} \\
				&\leq  \frac{3}{\epsilon}\norm{N_i^c}_{C^0(B_{\epsilon}(0))}  +\norm{DN_i^c}_{C^0(B_{\epsilon}(0))} .
			\end{align*} 
			\noindent Thus
			\begin{align*}
				\sup_w\norm{(DM)_w}_{op,0}& = \sup_w\max_i(\norm{(DM_i)_{w_i}}_{op}) \\
				&\leq  \max_i\left(\frac{3}{\epsilon}\norm{N_i^c}_{C^0(B_{\epsilon}(0))}  +\norm{DN_i^c}_{C^0(B_{\epsilon}(0))}\right) . 
			\end{align*} So \ref{claim:second} will follow from the following two subclaims. 
			\begin{claim}
				$\norm{N_i^c}_{C^0(B_{\epsilon}(0))}\leq \delta+2\epsilon^2K$.
			\end{claim}
			\noindent Well, for all $v\in B_{\epsilon}(0)$, Taylor's theorem implies that there exist $z_1, z_2\in B_{\epsilon}(0)$ such that
			\begin{align*}
				N_i^c(v)- h_i^c(0) = \frac 1 2
				\begin{pmatrix}
					v^t & 0 \\
					0 & v^t
				\end{pmatrix}
				\cdot
				\begin{pmatrix}
					(D^2 (h_i^c)_1)_{z_1} & 0 \\
					0 & (D^2 (h_i^c)_2)_{z_2}
				\end{pmatrix}
				\cdot \begin{pmatrix} v \\ v \end{pmatrix}
			\end{align*} where $h_i^c = ((h_i^c)_1, (h_i^c)_2)$. Thus \[\norm{ N_i^c(v)}_2\leq |v|^2\max_{k = 1,2}\norm{(D^2( h_i^c)_k)_{z_k}}_F + \delta.\] Since \[\max_{k = 1,2}\norm{(D^2( h_i^c)_k)_{z_k}}_F\leq 2 \norm{D^2 h_i^c}_{C^0(B_{\epsilon}(0)), \infty}\] and $|v| < \epsilon$, we obtain  
			\begin{align*}
				\norm{N_i^c}_{C^0(B_{\epsilon}(0))}&\leq\delta+2\epsilon^2K.
			\end{align*}
			\par\medskip
			\begin{claim}
				$\norm{DN_i^c}_{C^0(B_{\epsilon}(0)),op}\leq2\sqrt{2}\epsilon K$.
			\end{claim}
			\noindent Notice
			\begin{align*}
				(DN_i^c)_w  = (Dh_i^c)_w - (Dh_i^c)_{0}.
			\end{align*}
			By the mean value theorem applied to the four entries of $Dh_i^c$,
			\begin{align*} \label{idk}
				\norm{ (DN_i^c)_w}_{op} = \norm{(Dh_i^c)_w - (Dh_i^c)_{0}}_{op} \leq \norm{(D h_i^c)_w - (D h_i^c)_{0}}_F \leq2\sqrt{2}\epsilon K \notag
			\end{align*}
			Combining the two subclaims gives
			\begin{align*}
				\sup_{w\in\R^{2n}}\norm{(DM)_w}_{op,0}
				&\leq\frac{3}{\epsilon}\left(\delta+2\epsilon^2K\right)
				+2\sqrt{2}\epsilon K =\frac{3\delta}{\epsilon}+(6+2\sqrt{2})\epsilon K <\frac{3\delta}{\epsilon}+9\epsilon K,
			\end{align*}
			proving Claim \ref{claim:second}, and hence Claim \ref{claim:first}.
			We are not quite done, since we need $z\in B_{\frac{ \epsilon}{2}}(\vec 0)$. Well, $z = \lim_{i\to\infty}\mathcal{F}^i(\hat z)$ for any $\hat z\in \R^{2n}$. 
			Therefore,
			\begin{align*}
				\norm{z- \hat z}_0& \leq \sum_i \norm{\mathcal{F}^{i+1}(\hat z) - \mathcal{F}^i(\hat z)}_0  \\
				&\leq \sum_i \left(\frac{13}{16}\right)^i\norm{\mathcal{F}(\hat z) - \hat z}_0 =\frac{16}{3}\norm{\mathcal{F}(\hat z) - \hat z}_0 .
			\end{align*}
			In particular, substituting $\vec 0$ for $\hat z$, 
			\begin{align*}
				\norm{z}_0
				\leq\frac{16}{3}\norm{\mathcal{F}(0)}_0
				&\leq\frac{16}{3}\kappa\norm{M(0)}_0\\
				&=\frac{16}{3}\kappa\norm{H(0)-L(0)}_0\\
				&=\frac{16}{3}\kappa\norm{H(0)}_0 <\frac{16}{3}\kappa\delta
				<6\kappa\delta
				=\frac{\epsilon}{2},
			\end{align*}
			concluding the proof.
		\end{proof}
		\section{Miscellaneous}
		
		\begin{lem}\label{lem:smoothXA}
			For every nonzero $A\in\R$, the variety $X_A$ is smooth.
			\begin{proof}
				We first check that $X_A$ is smooth on its affine part. Suppose that $(x,y,z)$ were a singular point, meaning that $\partial_xq_A$, $\partial_yq_A$, and $\partial_zq_A$ vanish at $(x,y,z)$. First, notice $xyz\neq0$. 
				Indeed, if, say, $x=0$, then $\partial_xq_A=Ayz=0$. After interchanging $y$ and $z$ if necessary, we may assume that $y=0$, and then $\partial_zq_A=2z=0$. But $q_A(0,0,0)$ is nonzero, a contradiction. Moreover, none of $1+x^2$, $1+y^2$, and $1+z^2$ vanishes. Indeed, if $1+z^2=0$, then $\partial_xq_A=Ayz\neq0$. Now
				\begin{align*}
					x\partial_xq_A-y\partial_yq_A
					=2(x^2-y^2)(1+z^2),
				\end{align*}
				and the analogous identities show that
				\begin{align*}
					x^2=y^2=z^2=s.
				\end{align*}
				Now $x\partial_xq_A=0$ implies $Axyz=-2s(1+s)^2$. Substituting into $q_A=0$, and squaring $Axyz=-2s(1+s)^2$, gives
				\begin{align*}
					s^3+s^2-s+1=0, \hspace{1 cm} A^2s=4(1+s)^4.
				\end{align*}
				Eliminating $s$ yields
				\begin{align*}
					A^6-48A^4+896A^2+1024=0,
				\end{align*}
				which one can check has no real solutions. 
				
				Now assume that $x_0=0$, and set $x=1$. On $X_A$, we then have
				\begin{align*}
					(y_0^2+y^2)(z_0^2+z^2)=0.
				\end{align*}
				If, say, $y_0^2+y^2 = 0$ and $z_0^2+z^2\neq 0$, then $\partial_{y_0}\tilde q_A = 0$ implies $y_0 = 0$, and hence $y = 0$, a contradiction. If both factors vanish, then
				\begin{align*}
					\frac{\partial \tilde q_A}{\partial x_0}
					=A y_0yz_0z\neq0,
				\end{align*} a contradiction. Hence $X_A$ is smooth outside of the affine part. We conclude that $X_A$ is smooth.
				
			\end{proof}
		\end{lem}
		
		Recall from Notation \ref{notation} that we let $X\coloneqq X_{10}$, $q\coloneqq q_{10}$, and $\tilde q\coloneqq \tilde q_{10}$.
		\begin{lem}\label{lem:star} 
			$X(\R)$ is completely contained in the affine part of $(\P^1)^3$. Moreover, if $S^2\subset \R^3$ denotes the unit sphere, then radial projection $R:X(\R)\to S^2$ is a diffeomorphism.
			\begin{proof}
				In homogeneous coordinates, we find
				\[\tilde q([0:1], [y_0: y], [z_0: z]) = (y_0^2 + y^2)(z_0^2 + z^2),\]
				so, by the symmetry of $\tilde q$, $X(\R)$ contains no points at infinity. 
				
				If every radial intersection, that is, an intersection of a ray from the origin with $X(\R)$, were transverse, then $R:X(\R)\to S^2$ would be a local diffeomorphism. In that case, the image of $R$ would be a nonempty subset of $S^2$ that is both closed and open, so $R$ would be surjective. Moreover, since $X(\R)$ is compact, $R$ would be a covering map. Now the positive $x$-axis meets $X(\R)$ exactly once, at $(1,0,0)$, since $q(t,0,0)=t^2-1$. Thus this covering would be one-sheeted and $R$ would be a diffeomorphism.

				For the sake of contradiction, assume there exists $p\in X(\R)$ at which the radial intersection is tangent, or equivalently $p\in T_pX(\R)$. In other words, $q(p)=0$ and $Dq_p(p)=0$, where
				\[	Dq_p(p)=x\frac{\partial q}{\partial x}(p)
					+y\frac{\partial q}{\partial y}(p)
					+z\frac{\partial q}{\partial z}(p).\]
				Write $p=(x,y,z)$ and set $a=x^2$, $b=y^2$, $c=z^2$, $t=xyz$, and 
				\[ s_1=a+b+c, \hspace{1 cm} s_2=ab+bc+ca.\]
				Then $q(p)=0$ gives
				\begin{align} \label{eq:radial-q}
					s_1+s_2+t^2+10t=1,
				\end{align}
				and $Dq_p(p)=0$ gives
				\begin{align}\label{eq:radial-Dq}
					2s_1+4s_2+6t^2+30t=0. 
				\end{align}
				Eliminating the term linear in $t$ from \eqref{eq:radial-q} and \eqref{eq:radial-Dq} yields
				\[s_1=3+s_2+3t^2,\]
				so
				\begin{align}\label{thingy}
					10t=-2(1+s_2+2t^2).
				\end{align}
Then \eqref{thingy} becomes
	\begin{align*}
		-5t=1+s_2+2t^2.
	\end{align*}  It remains to eliminate $s_2$.
Now, since $ab, bc, ca\geq 0$, we can compare their arithmetic and geometric means to obtain
	\begin{align}\label{eq:radial-s2}
	s_2=ab+bc+ca \geq 3((ab)(bc)(ca))^{1/3} =3t^{4/3}. 
\end{align} Hence
	\begin{align*}
	-5t\geq 1+3t^{4/3} +2t^2.
\end{align*} Letting $u = t^{1/3}$, we have 
	\begin{align*}
	1+ 5u^3 + 3u^4 +2u^6  \leq 0.
\end{align*} But the above polynomial is positive for all $u$, a contradiction.
			\end{proof}
		\end{lem} 
		
		\begin{lem}\label{claim:1} We have
			\begin{align*}
				\max_{(x,y,z)\in X(\R)}|x|<2.31.
			\end{align*}
			By symmetry, the same bound holds for the $y$ and $z$ coordinates.
			\begin{proof}
				Since $(x,y,z)\in X(\R)$, we have $x,y\in \R$, and the discriminant of $q$ as a function of $z$ is non-negative. Hence
				\begin{align*}
					100x^2y^2-4(1+x^2)(1+y^2)((1+x^2)(1+y^2)-2)\geq0.
				\end{align*}
				Setting $s=1+y^2$ and $u  = x^2$ yields
				\begin{align*}
					(1+u)^2s^2-(27u+2)s+25u\leq0.
				\end{align*}
				Now this is a quadratic polynomial in $s$, nonpositive at $s=1+y^2$, with positive leading coefficient. Hence it has a real zero and nonnegative discriminant. So
				\begin{align*}
					(27u+2)^2-100u(1+u)^2\geq0,
				\end{align*}
				which expands to
				\begin{align*}
					100u^3-529u^2-8u-4\leq0.
				\end{align*}
				Now this polynomial in $u$ has positive leading coefficient, and its largest zero is approximately $5.306$. Hence $u < 5.31$, and $|x| < \sqrt{5.31} < 2.31$, concluding the lemma.
			\end{proof}
		\end{lem}
		
		\begin{lem}\label{lem:2} 
			Let $S = \set{(x,y)\in \R^2\mid \mathcal{D}(x,y)>0}$.
			The function $\mathcal{D}$ satisfies the derivative bounds 
			\[\max_{|\beta| = 3} \sup_{(x,y)\in S} |D^\beta \mathcal{D}(x,y)| \leq 20000.\]
			\begin{proof}
				We find 
				\[\partial_{xxx}\mathcal{D}(x,y) = -96 x (1 + y^2)^2\] and 
				\begin{align*}
					\partial_{xxy}\mathcal{D}(x,y) =432 y - 128 x^2 y (1 + y^2) - 64 (1 + x^2) y (1 + y^2).
				\end{align*}
				If $|x|,|y|\leq c$, we calculate
				\[|\partial_{xxx}\mathcal{D}(x,y)| \leq 96 c (1 + c^2)^2\] and 
				\begin{align*}
				|\partial_{xxy}\mathcal{D}(x,y)| \leq432 c + 128 c^3 (1 + c^2) + 64  c(1 + c^2)^2. 
				\end{align*}
				The bound from Lemma \ref{claim:1} and symmetry conclude the proof.
			\end{proof}
		\end{lem} 
		
		\begin{lem}\label{lem:well}
			Considering $f_A$ as a smooth function from an open neighborhood of $X_A(\R)$ into $\R^3$, we have
			\[\norm{D^2f_A}_{C^0(X_A(\R)), \infty} \leq 9^2(2\cdot \max(|A|, 1)^2)^3\]
			and 
			\[\norm{Df_A}_{C^0(X_A(\R)),\infty} \leq 3^2\cdot \max(|A|, 1)^3.\]
			
			\begin{proof} Define 
				\begin{align*}
					\alpha(x,y) \coloneqq \frac{xy}{(1+x^2)(1+y^2)}.
				\end{align*} In affine coordinates, one calculates
			\begin{align*}
				\sigma_1^A(x,y,z) = \left(-x - A\alpha(y,z), y, z\right), \hspace{1cm} \sigma_2^A(x,y,z) = \left(x, -y-A\alpha(x,z), z\right),
			\end{align*}  and 
			\[\sigma_3^A(x,y,z) =\left(x,y, -z - A\alpha(x,y)\right).\]

			 Furthermore, notice that 
				\begin{align*}
					\alpha(x,y) = \frac{1}{4}\left(\frac{1}{-x - i} + \frac{1}{-x + i}\right)\left(\frac{1}{-y - i} + \frac{1}{-y + i}\right).
				\end{align*}
				It follows that, for $x,y\in \R$,
				\begin{align*}
					|\alpha| \leq \frac 1 4 (1 + 1)(1 + 1) = 1, \hspace{1 cm} |\partial_y(\alpha)|, |\partial_x(\alpha)| \leq \frac 1 4 (1 + 1)(1 + 1) = 1
				\end{align*}
			and 
			\begin{align*}
						|\partial_{xy}(\alpha)|\leq \frac 1 4 (1 + 1)(1 + 1) = 1,  \hspace{1 cm} |\partial_{yy}(\alpha)|,  |\partial_{xx}(\alpha)|\leq \frac 1 4 2(1 + 1)(1 + 1) = 2.
				\end{align*}
			
				Then $\norm{D \sigma_i^A}_{X_A(\R), \infty} \leq \max(1, |A|)$, and, by \eqref{inequality}, $\norm{J^2(\sigma_i^A)}_{X_A(\R), \infty} \leq 2\cdot \max(|A|, 1)^2$.
				The result follows from matrix multiplication and dimension counting.
			\end{proof}
			
		\end{lem}
		\begin{lem} \label{lem:7}
			If $C_1, C_2\in \R_{>0}$ and $x,y\in \R$ are such that 
			\begin{itemize}
				\item $1 \leq \mathcal{D}(x,y) \leq C_1$ and 
				\item $\norm{(D  \mathcal{D})_{(x,y)}}_{\infty}, \norm{(D^2  \mathcal{D})_{(x,y)}}_{\infty} \leq C_2$
			\end{itemize} then
			\begin{itemize}
				\item $\norm{(D^2 p_\pm)_{(x,y)}}_{\infty}\leq 10 +C_1^{0.5} +  \frac 3 4 C_2 + \frac 1 8 C_2^2$ and 
				\item $\norm{(Dp_\pm)_{(x,y)}}_{\infty} \leq 5 + \frac 1 2 C_1^{0.5} + \frac 1 4 C_2.$
			\end{itemize}
			\begin{proof} With $\alpha$ as in Lemma \ref{lem:well} and $\beta(x,y)\coloneqq \frac{1}{(1+x^2)(1+y^2)}$,
				\begin{align*}
					p_\pm & =-5\alpha \pm \frac 1 2 \mathcal{D}^{\frac{1}{2}} \beta.
				\end{align*} Notice that the same derivative bounds we established for $\alpha$ in Lemma \ref{lem:well} hold for $\beta$. Then
				\begin{align*}
					\left|\partial_j p_\pm\right| & \leq 5\left|\partial_j \alpha\right| +  \frac 1 2 \left(\left|\mathcal{D}^{0.5}\right| \left|\partial_j \beta\right| + \left|\beta\right| \left|\frac{\partial_j \mathcal{D}}{2 \mathcal{D}^{0.5}}\right| \right) \\ 
					& \leq 5 + \frac 1 2  C_1^{0.5} + \frac 1 4 C_2
				\end{align*} and
				\begin{align*}
					\left|\partial_{ij} p_\pm\right| &\leq 5\left|\partial_{ij} \alpha\right| +  \frac 1 2 \left(\left|\frac{\partial_i \mathcal{D}}{2\mathcal{D}^{0.5}}\right|\left|\partial_j \beta\right| + \left|\mathcal{D}^{0.5}\right|\left|\partial_{ij}\beta\right| + \left|\beta\right| \left|\frac{\partial_j \mathcal{D} \partial_i \mathcal{D}}{4 \mathcal{D}^{1.5}}\right| + \left|\beta\right| \left|\frac{\partial_{ij}\mathcal{D}}{2 \mathcal{D}^{0.5}}\right| + \left|\partial_i \beta\right| \left|\frac{\partial_j \mathcal{D}}{2\mathcal{D}^{0.5}}\right| \right) \\
					& \leq 10 +  \frac 1 2 \left(\frac{C_2}{2} + 2 C_1^{0.5} + \frac{C_2^2}{4} + \frac{C_2}{2} + \frac{C_2}{2} \right) \\
					& = 10 + C_1^{0.5} +  \frac 3 4 C_2 + \frac 1 8 C_2^2,
				\end{align*} as desired. 
			\end{proof}
		\end{lem}
		
		\subsection{Jets}\label{sec:jets}
		Consider $C^\infty(\R^n,\R)$, the collection of smooth functions $f:\R^n\to \R$. Let $k\in \N$ and $x_0\in \R^n$, and consider the equivalence relation on $C^\infty(\R^n,\R)$ defined by $f\sim g$ if $f(x_0) = g(x_0)$ and all partials of $f$ and $g$ at $x_0$ of order $\leq k$ agree. We call such an equivalence class a $k$-jet, and denote the set of classes by 
		$J^k_{x_0}(\R^n, \R)$. Furthermore, denote by $J^k_{x_0}(\R^n, \R)_0$ the collection of jets mapping $x_0$ to $0$. By Taylor's theorem, $J^k_{x_0}(\R^n, \R)_0$ is a linear space.

		For a multi-index $\alpha=(\alpha_1,\ldots,\alpha_n)\in\Z_{\geq0}^n$, write
		\begin{align*}
			|\alpha|\coloneqq\sum_{j=1}^n\alpha_j, \hspace{1 cm}
			\alpha!\coloneqq\prod_{j=1}^n\alpha_j!, \hspace{1 cm}
			(x-x_0)^\alpha\coloneqq
			\prod_{j=1}^n\big(x_j-(x_0)_j\big)^{\alpha_j}.
		\end{align*}
		The normalized monomials
		\begin{align*}
			\mathcal{B}_{x_0}
			\coloneqq
			\set{
				\frac{(x-x_0)^\alpha}{\alpha!}
				\ \middle|\
				\alpha\in\Z_{\geq0}^n,\ 1\leq|\alpha|\leq k
			}
		\end{align*}
		form a basis for $J^k_{x_0}(\R^n, \R)_0$. We use $J^k(\R^n, \R)_0$ to denote the trivial bundle obtained by assigning the vector space $J^k_{x_0}(\R^n, \R)_0$ to each $x_0\in \R^n$.
		
		For any $f\in C^k(\R^n, \R^m)$ and $x_0\in \R^n$, $f$ induces a linear map $J^k_{x_0}(f):
			J_{f(x_0)}^k(\R^m, \R)_0
			\longrightarrow
			J_{x_0}^k(\R^n, \R)_0$ by precomposition.
		The entries of $J^k_{x_0}(f)$ depend on the partials of $f$ at $x_0$ up to order $k$. Moreover, jets behave nicely under composition; precisely, whenever $f\circ g$ is defined, $J^k_{x_0}(f\circ g)
			=
			J^k_{x_0}(g)\cdot J^k_{g(x_0)}(f)$ in the sense of matrix multiplication.
		Often, we use $J^k(f)$ to denote the matrix-valued function $x_0\mapsto J^k_{x_0}(f)$.
		For example, consider a smooth map $f:\R^2\to \R^2$ with coordinate functions $f_1, f_2$. Order the basis $\mathcal B_{x_0}$ as
		\begin{align*}
			\left(
			x_1-(x_0)_1,\,
			x_2-(x_0)_2,\,
			\frac{(x_1-(x_0)_1)^2}{2},\,
			\frac{(x_2-(x_0)_2)^2}{2},\,
			(x_1-(x_0)_1)(x_2-(x_0)_2)
			\right),
		\end{align*}
		and use the analogous ordering at $f(x_0)$. The $2$-jet $J^2(f)$, expressed with respect to the bases $\mathcal B_{x_0}$ and $\mathcal B_{f(x_0)}$ on $J^2_{x_0}(\R^2,\R)_0$ and $J^2_{f(x_0)}(\R^2,\R)_0$, respectively, is
		\begin{align*} J^2(f) = 
			\begin{pmatrix}
				\frac{\partial f_1}{\partial x_1} & \frac{\partial f_2}{\partial x_1}  & 0 & 0 & 0 \\
				\frac{\partial f_1}{\partial x_2} & \frac{\partial f_2}{\partial x_2}  & 0 & 0 & 0 \\
				\frac{\partial^2 f_1}{\partial x_1^2} & \frac{\partial^2 f_2}{\partial x_1^2}  & \left( \frac{\partial f_1}{\partial x_1} \right)^2 & \left( \frac{\partial f_2}{\partial x_1} \right)^2 & 2\frac{\partial f_1}{\partial x_1}  \frac{\partial f_2}{\partial x_1} \\
				\frac{\partial^2 f_1}{\partial x_2^2} & \frac{\partial^2 f_2}{\partial x_2^2} & \left( \frac{\partial f_1}{\partial x_2} \right)^2 & \left( \frac{\partial f_2}{\partial x_2} \right)^2 & 2\frac{\partial f_1}{\partial x_2} \frac{\partial f_2}{\partial x_2}  \\
				\frac{\partial^2 f_1}{\partial x_1 \partial x_2} & \frac{\partial^2 f_2}{\partial x_1 \partial x_2} &  \frac{\partial f_1}{\partial x_1} \frac{\partial f_1}{\partial x_2} & \frac{\partial f_2}{\partial x_1} \frac{\partial f_2}{\partial x_2} & \frac{\partial f_1}{\partial x_1} \frac{\partial f_2}{\partial x_2} + \frac{\partial f_2}{\partial x_1} \frac{\partial f_1}{\partial x_2}
			\end{pmatrix}.
		\end{align*} 
		In general, one has 
		\begin{align}\label{inequality}
			\norm{J^2_{x_0}(f)}_{\infty}
			\leq
			\max\left(
				\norm{(Df)_{x_0}}_{\infty},
				\norm{(D^2f)_{x_0}}_{\infty},
				2\norm{(Df)_{x_0}}^2_{\infty}
			\right).
		\end{align}

		\section{Computer-assisted certifications}\label{sec:errors}

		\subsection{MPFR interval arithmetic}
		The following interval-arithmetic framework, based on the GNU MPFR library \cite{MR2326955}, applies to \texttt{verifyPeriodicOrbit.c}, \texttt{derivativesInterval.c}, and \texttt{verifyArcs.c}. A nonzero finite MPFR number of precision $N$ has the form
		\begin{align*}
			\pm 2^{p}\left(\frac{1}{2}+\sum_{k=2}^{N}c_k2^{-k}\right),
		\end{align*}
		where $c_k\in\set{0,1}$ and $p$ is an integer. The library allows users to specify an exponent range $[m,M]$. For fixed $N$ and $[m,M]$, we call zero and the numbers of the displayed form with $p\in[m,M]$ \emph{representable}. The smallest positive representable number is $2^m/2$, and the largest is $2^M(1-2^{-N})$. We say that a real number $x$ is in the \emph{allowable range} if $x=0$ or
		\begin{align*}
			2^m/2\leq |x|\leq 2^M(1-2^{-N}).
		\end{align*}

		For a real number $x$ in the allowable range, write $\operatorname{rd}(x)$ and $\operatorname{ru}(x)$ for the representable values obtained by rounding $x$ toward $-\infty$ and $+\infty$, respectively. Each elementary MPFR operation takes allowable inputs and returns the exact real result rounded to a representable value according to the chosen mode, provided that the exact result lies in the allowable range. The modes \texttt{MPFR\_RNDD} and \texttt{MPFR\_RNDU} correspond to downward and upward rounding, respectively.

		In this framework, a real number is represented by an interval containing it whose endpoints are representable. The programs then perform \emph{outwardly rounded} interval arithmetic, which is standard interval arithmetic in which lower endpoints are rounded downward and upper endpoints are rounded upward. This ensures that every computed interval contains the corresponding exact result. Precisely, for intervals $I=[a,b]$ and $J=[c,d]$, the programs use the following:
		\begin{align*}
			&I+J
			\subset
			[\operatorname{rd}(a+c),\operatorname{ru}(b+d)], \hspace{1 cm}
			I-J \subset
			[\operatorname{rd}(a-d),\operatorname{ru}(b-c)], 
		\end{align*}
		and 
		\begin{align*}
			IJ \subset
			[
			\min_{s\in\set{a,b}, \ t\in\set{c,d}}
			\operatorname{rd}(st),
			\max_{s\in\set{a,b}, \ t\in\set{c,d}}
			\operatorname{ru}(st)
			].
		\end{align*}
		If $0\notin J$, it computes
		\begin{align*}
			J^{-1} \subset
			[\operatorname{rd}(1/d),\operatorname{ru}(1/c)], \hspace{1 cm}
			I/J=IJ^{-1}.
		\end{align*}
		For $0<a\leq b$, it computes
		\begin{align*}
			\sqrt I
			\subset
			[\operatorname{rd}(\sqrt a),\operatorname{ru}(\sqrt b)].
		\end{align*}
	Before performing interval division, the programs require the divisor interval to exclude zero. Before evaluating a square root, they require the interval under the radical to be nonnegative. If either requirement cannot be verified, the relevant program terminates.

		\subsection{Certification in Sections \ref{sec:1} and \ref{sec:kappa} --  \texttt{derivativesInterval.c}}\mbox{}\par
		The program \texttt{derivativesInterval.c} uses the preceding interval operations with precision $N=5000$. The program uses the polynomial formulas for $\mathcal D$ and its first and second partials to compute certified enclosures at each $(a_i,b_i)$. These enclosures certify that the corresponding values displayed in Table \ref{table:3} are correctly rounded to two decimal places. Then Table \ref{table:3} gives us the bounds $K<115$, $R<163$ and $M<442$. 
		
		The same evaluation on the square of radius $10^{-3}$ about each $(a_i,b_i)$ certifies that $\mathcal D>0$ there, and hence on $B_{10^{-3}}(a_i,b_i)$.

		To estimate $L_i=(Df_i^c)_0$, the program uses explicit derivative formulas for the charts and the three involutions, combined with interval matrix multiplication, to certify enclosures for each entry of $L_i$. If $I_{i,jk}$ is the resulting enclosure of the $(j,k)$-entry of $L_i$, the program certifies 
		\begin{align*}
			\max_{i,j,k}\sup_{x\in I_{i,jk}}|x-(\widetilde L_i)_{jk}| < 10^{-2}.
		\end{align*}
	Hence $\norm{L_i-\widetilde L_i}_\infty < 10^{-2}$, which is what we used in Section \ref{sec:kappa}.
		\subsection{Certifying the pseudo-orbit -- \texttt{verifyPeriodicOrbit.c}}
		The program \texttt{verifyPeriodicOrbit.c} uses the preceding interval operations with precision $N=500$. The remaining certifications needed in Section \ref{sec:ApplyShadowing} are that $0\in U_i$ for every $0\leq i\leq9$ and that
			\begin{align}\label{eq:certifiedResidualGoal}
			\max_{0\leq i\leq9}\norm{f_i^c(0)}_2<10^{-29}.
		\end{align}  Recall that $\xi_i:\R^3\to \R^2$ is the affine coordinate projection that inverts $\phi_i$ on $\phi_i(V_i)$.
			For the remainder of this subsection, indices are taken modulo $10$. 
			Let $\tilde\phi_i=\Psi_{j_i}^{\epsilon_i}$ be the chart listed in Table \ref{table:1}. The program certifies a box containing $f\big(\tilde\phi_i(a_i,b_i)\big)$ using the formulas for $\mathcal D$, $p_\pm$, and $\Psi_j^\pm$ from Section \ref{sec:charts}, and then applies $f=\sigma_3\circ\sigma_2\circ\sigma_1$.
			
				Let $(u,v)\coloneqq\xi_{i+1}\big(f(\tilde\phi_i(a_i,b_i))\big)+(a_{i+1},b_{i+1})$ and let $w$ be the coordinate lost in projection. The program first certifies that $\mathcal D(u,v)>0$. It also checks that the image lies on the branch of $X(\R)$ parametrized by $\tilde\phi_{i+1}$. The branch is determined by the sign of
				\begin{align*}
					2(1+u^2)(1+v^2)w+10uv.
				\end{align*}
			The program certifies that its interval enclosure has the sign corresponding to $\tilde\phi_{i+1}$. Thus $f(\phi_i(0))\in\phi_{i+1}(V_{i+1})$, and hence $0\in U_i$ for every $i$.
			
			We now have
				\begin{align}\label{eq:residualEvaluation}
					f_i^c(0) = \xi_{i+1}(f(\tilde\phi_i(a_i,b_i)))
		\end{align} and intervals $R_{i,1}$ and $R_{i,2}$ satisfying $f_i^c(0)\in R_{i,1}\times R_{i,2}$.
		
		For an interval $R=[\ell,u]$, define
		\begin{align*}
			m(R)\coloneqq\max\set{|\ell|,|u|}.
		\end{align*}
		Using upward rounding, the program then computes
		\begin{align*}
			N_i
			\coloneqq
			\operatorname{ru}
			\left(
			\sqrt{m(R_{i,1})^2+m(R_{i,2})^2}
			\right),
		\end{align*}
		which bounds $\norm{f_i^c(0)}_2$ from above.
		
		Now Table \ref{table:residualCertificate} records the output of \texttt{verifyPeriodicOrbit.c}. Here $\mathcal D_i^{\mathrm{lb}}$ is a lower bound for the discriminant used to construct $\tilde\phi_i(a_i,b_i)$, while $M_{i,k}$ and $N_i$ satisfy
		\begin{align*}
			R_{i,k}\subset[-M_{i,k},M_{i,k}] \hspace{.5 cm} \text{and}\hspace{.5 cm} \norm{f_i^c(0)}_2\leq N_i.
		\end{align*}
			In Table \ref{table:residualCertificate}, the displayed values of $\mathcal D_i^{\mathrm{lb}}$ have been rounded downward, while all other bounds have been rounded upward.
	\begin{table}[h]
	\centering
	\begin{tabular}{c|cccc}
	$i$ & $\mathcal D_i^{\mathrm{lb}}$ & $M_{i,1}$ & $M_{i,2}$ & $N_i$\\
	\hline
	$0$ & $114$  & $5.67\cdot10^{-32}$ & $3.51\cdot10^{-31}$ & $3.56\cdot10^{-31}$\\
	$1$ & $8.70$ & $3.17\cdot10^{-31}$ & $1.00\cdot10^{-30}$ & $1.05\cdot10^{-30}$\\
	$2$ & $70.2$ & $1.29\cdot10^{-30}$ & $7.40\cdot10^{-31}$ & $1.49\cdot10^{-30}$\\
	$3$ & $8.59$ & $1.09\cdot10^{-30}$ & $1.67\cdot10^{-31}$ & $1.10\cdot10^{-30}$\\
	$4$ & $108$  & $1.67\cdot10^{-31}$ & $9.62\cdot10^{-32}$ & $1.93\cdot10^{-31}$\\
	$5$ & $10.4$ & $1.90\cdot10^{-31}$ & $4.12\cdot10^{-31}$ & $4.53\cdot10^{-31}$\\
	$6$ & $108$  & $1.76\cdot10^{-31}$ & $3.69\cdot10^{-31}$ & $4.09\cdot10^{-31}$\\
	$7$ & $8.59$ & $7.92\cdot10^{-31}$ & $7.05\cdot10^{-30}$ & $7.09\cdot10^{-30}$\\
	$8$ & $70.2$ & $2.00\cdot10^{-30}$ & $2.58\cdot10^{-30}$ & $3.26\cdot10^{-30}$\\
	$9$ & $8.70$ & $5.73\cdot10^{-31}$ & $7.46\cdot10^{-31}$ & $9.40\cdot10^{-31}$
\end{tabular}
\vspace{.25 cm}\\
	\caption{Certified interval bounds for the pseudo-orbit residuals.}
\label{table:residualCertificate}
\end{table}
		In particular,
		\begin{align*}
			\max_{0\leq i\leq9}\norm{f_i^c(0)}_2
			\leq
			\max_{0\leq i\leq9}N_i
			\leq7.09\cdot10^{-30}
			<10^{-29},
		\end{align*}
		which proves \eqref{eq:certifiedResidualGoal}.

		\subsection{Certifying relative homotopy classes using arc-data -- \texttt{verifyArcs.c}}\label{sec:arc-certification}
				It remains, for $0\leq i\leq8$, to certify that $f'\circ p_i$ has the relative homotopy class represented by the corresponding arc in the top panel of Figure \ref{fig:actual}. In this section, we describe a strategy to do so, which is implemented in \texttt{verifyArcs.c}.
				The program \texttt{verifyArcs.c} uses the preceding interval operations with precision $N=500$.
			
					Recall from Section \ref{sec:mclass} that $\psi_\gamma:S^2\setminus\{\gamma\} \to \R^2$ denotes stereographic projection through $\gamma$. Throughout this subsection, we work in the stereographic coordinate plane $\R^2$, identified topologically with $\mathring D$. Let $q_\gamma\coloneqq(\psi_\gamma\circ R\circ f^{-1}\circ R^{-1})(\gamma)$ and define $\Pi:\R^2\setminus\set{q_\gamma} \to \R^2$ by 
				\begin{align*}
					\Pi\coloneqq \psi_\gamma \circ R \circ f\circ R^{-1}\circ \psi_\gamma^{-1},
				\end{align*}
				where $R:X(\R)\to S^2$ is the restriction to $X(\R)$ of radial projection. Recall that $f'$ is the homeomorphism of $\R^2$ obtained by stereographically conjugating a representative of $[f]$ on $S^2$ that fixes $\gamma$ and whose restriction to each lifted arc $\hat p_i=\psi_\gamma^{-1}\circ p_i$ agrees with that of $f$. Since the arcs $p_i$ avoid $q_\gamma$, it follows that $f'\circ p_i=\Pi\circ p_i$ for every $0\leq i\leq8$. Hence our program \texttt{verifyArcs.c} will work with the explicit map $\Pi$.

				 Recall from Section \ref{sec:mclass} that $z_0,\dots,z_9$ are the images under $\psi_\gamma\circ R$ of the points $f^j(x)$, indexed from left to right in the plane, and that $\tau$ is characterized by
				 \begin{align*}
				 	z_{\tau(j)}=(\psi_\gamma\circ R)(f^j(x)).
				 \end{align*}
			 The certification in \texttt{verifyArcs.c} happens in two steps:
			 
			 {\bf Step $1$.}
					 \texttt{verifyArcs.c} outputs the permutation $\tau$ and pairwise disjoint rational boxes $B_0,\dots,B_9\subset\R^2$ with $z_i\in\operatorname{int}(B_i)$. The output boxes are certified to satisfy
			 \begin{align*}
			 	\psi_\gamma\circ R\circ\phi_j([-\bar\eta,\bar\eta]^2)\subset\operatorname{int}(B_{\tau(j)})
			 \end{align*}
				 for all $0\leq j\leq9$, where $\bar\eta\coloneqq5\cdot10^{-27}$ is as in Section \ref{sec:ApplyShadowing}.
				 The program also certifies that every point of $B_i$ lies strictly to the left of every point of $B_{i+1}$ for $0\leq i\leq8$, and that $\gamma\notin R\circ\phi_j([-\bar\eta,\bar\eta]^2)$ for every $0\leq j\leq9$. Consequently, the points $z_i$ are pairwise distinct and $x$ has least period ten. 
				
			 {\bf Step $2$.} Fix $0\leq k\leq8$. Let $k_i,k_f\in\set{0,\dots,9}$ be the indices determined by $\Pi(z_k)=z_{k_i}$ and $\Pi(z_{k+1})=z_{k_f}$. Then \texttt{verifyArcs.c} finds an integer $m\geq2$, rational numbers
			 \begin{align*}
			 	0=t_0<t_1<\cdots<t_m=1,
			 \end{align*}
			 and rectangles $R_0,\dots,R_{m-1}$ which satisfy the following properties, each of which is certified by \texttt{verifyArcs.c}. For $0\leq r<m$, set $I_r=[t_r,t_{r+1}]$. Then
		\begin{enumerate}
			\item each $R_r$ has rational vertices;
			\item for every $0\leq r<m$ and every $t\in I_r$, the map $\Pi$ is defined on $(1-t)B_k+tB_{k+1}$ and
			\begin{align*}
				\Pi\big((1-t)B_k+tB_{k+1}\big)\subset R_r.
			\end{align*}
			\item for every $0\leq r<m$ and $0\leq\ell\leq9$,
			\begin{align*}
				R_r\cap B_\ell=\varnothing
				\end{align*}
				unless $(r,\ell)=(0,k_i)$ or $(r,\ell)=(m-1,k_f)$;
				\item $B_{k_i}\subset R_0$ and $B_{k_f}\subset R_{m-1}$.
			\end{enumerate}

			By Lemma \ref{lem:star}, the radial inverse $R^{-1}:S^2\to X(\R)$ sends each $u\in S^2$ to the unique point $tu\in X(\R)$ with $t>0$. On each input box $U\subset\R^2$, to compute a box containing $\Pi(U)$, \texttt{verifyArcs.c} first certifies an interval $T\subset(0,\infty)$ containing, for every $v\in\psi_\gamma^{-1}(U)$, the unique positive solution $t$ of $q(tv)=0$. Multiplying a box containing $\psi_\gamma^{-1}(U)$ by $T$ gives a box containing $R^{-1}(\psi_\gamma^{-1}(U))$. The program then applies $f$, $R$, and $\psi_\gamma$ successively to this box using outward-rounded interval arithmetic to obtain a box containing $\Pi(U)$.
		
		 For each $i$, choose a rational reference point $\bar z_i\in\operatorname{int}(B_i)\cap\Q^2$. Let $\theta:\R^2\to \R^2$ be a homeomorphism supported on $\cup_i B_i$ such that $\theta(z_i) = \bar z_i$. 
				\begin{lem}\label{lem:homotope}
	 			Assuming the existence of such $I_r$ and $R_r$, we prove that $f'\circ p_k=\Pi\circ p_k$ is homotopic in $\R^2$ relative to $S=\set{z_0,z_1,\dots,z_9}$ to $\theta^{-1}\circ L_k$ for some polygonal path $L_k$ with vertices in $\Q^2$.
		\end{lem}
		Before proving Lemma \ref{lem:homotope}, we explain how it is used to certify the relative homotopy class of $f'\circ p_k$. Since $\theta:(\R^2,S)\to(\R^2,\theta(S))$ is a homeomorphism of pairs, $\theta^{-1}$ identifies the homotopy class of $L_k$ relative to $\theta(S)$ with that of $\theta^{-1}\circ L_k$ relative to $S$. By Lemma \ref{lem:homotope}, the latter is the relative homotopy class of $f'\circ p_k$. Thus it suffices to determine the relative homotopy class of $L_k$. The program does this by computing its arc-data with respect to $\theta(S)=\set{\bar z_0,\bar z_1,\dots,\bar z_9}$ exactly. This certified arc-data is used to illustrate the relative homotopy class of $f'\circ p_k$ in the top panel of Figure \ref{fig:actual}.
		
		\begin{proof}[Proof of Lemma \ref{lem:homotope}]
					For each $0 \leq r < m-1$, note that $R_r\cap R_{r+1}$ is a non-empty rational rectangle since it contains $\Pi\big((1-t_{r+1})B_k+t_{r+1}B_{k+1}\big)$. The program outputs a point $v_r\in R_r\cap R_{r+1}\cap\Q^2$ for every $0\leq r<m-1$. We let $v_{-1}\coloneqq \theta(\Pi(p_k(0)))$ and $v_{m-1}\coloneqq \theta(\Pi(p_k(1)))$. Define $L_k$ on $I_r$, for $0\leq r<m$, to be the straight-line path from $v_{r-1}$ to $v_r$. The points $v_r$ are chosen so that $L_k$ is embedded, as is also verified using exact rational arithmetic. The vertices $v_r$ are also chosen so that $L_k$ satisfies the hypotheses of Definition \ref{def:arc-data}, as verified using exact rational arithmetic. Both $v_{r-1}$ and $v_r$ are contained in $R_r$, so by convexity, $L_k(I_r)\subset R_r$. Since $\theta$ is supported on $\bigcup_{\ell=0}^9B_\ell$, properties (3) and (4) give $\theta(R_r)\subset R_r$. Thus $(\theta\circ f'\circ p_k)(I_r)\subset R_r$ for $0\leq r<m$.

					Using straight-line homotopies supported in the puncture-free intersections $R_r\cap R_{r+1}$, first homotope $\theta\circ f'\circ p_k$ to a path $\alpha$ such that $\alpha(t_{r+1})=v_r$. For $1\leq r<m-1$, the rectangle $R_r$ is puncture-free, so convexity of $R_r$ gives a straight-line homotopy between $\alpha|_{I_r}$ and $L_k|_{I_r}$ relative to $v_{r-1}$ and $v_r$. For $r=0$ and $r=m-1$, the rectangle $R_r$ contains no marked point other than the common marked endpoint of the two paths. Therefore, for $r\in\set{0,m-1}$, the paths $\alpha|_{I_r}$ and $L_k|_{I_r}$ are homotopic relative to $v_{r-1}$ and $v_r$. Combining these homotopies gives a homotopy relative to $\theta(S)$ from $\theta\circ f'\circ p_k$ to $L_k$. Composing it with $\theta^{-1}$ gives a homotopy in $\R^2$ relative to $S$ from $f'\circ p_k$ to $\theta^{-1}\circ L_k$, concluding the proof.
		\end{proof}

		Finally, \texttt{verifyArcs.c} use exact rational arithmetic to compute the arc-data of each polygonal path $L_k$ with respect to $\theta(S)=\set{\bar z_0,\bar z_1,\dots,\bar z_9}$, performs homotopy-preserving simplifications of this arc-data, and verifies that both the original and simplified sequences agree with the corresponding rows stored in \texttt{arc\_certificate\_data.c}. The program \texttt{mclass.c} uses the simplified sequences.

			\section{Tables}
		\begin{table}[h]
			\centering
			\begin{tabular}{c|cccccc}
				$i$ & $\mathcal{D}(a_i, b_i)$ & $\partial_1 \mathcal{D}(a_i, b_i)$ & $\partial_2 \mathcal{D}(a_i, b_i)$ & $\partial_{11} \mathcal{D}(a_i, b_i)$ & $\partial_{22} \mathcal{D}(a_i, b_i)$ & $\partial_{12} \mathcal{D}(a_i, b_i)$ \\
				\hline
				$0$ & 114.24 & 136.59 & -45.98 & -419.46 & -441.48 & -178.78 \\
				$1$ & 8.71 & -22.00 & 15.52 & 39.45 & 13.84 & -81.14 \\
				$2$ & 70.23 & 86.57 & -124.38 & -58.60 & -2.96 & -171.09 \\
				$3$ & 8.59 & -40.85 & 4.26 & 113.57 & -24.69 & -87.82 \\
				$4$ & 108.31 & 39.04 & 162.94 & -260.04 & -171.83 & 27.57 \\
				$5$ & 10.43 & 23.28 & -23.28 & 28.55 & 28.55 & -92.99 \\
				$6$ & 108.31 & -162.94 & -39.04 & -171.83 & -260.04 & 27.57 \\
				$7$ & 8.59 & -4.26 & -40.85 & -24.69 & 113.57 & 87.82 \\
				$8$ & 70.23 & 124.38 & -86.57 & -2.96 & -58.60 & -171.09 \\
				$9$ & 8.71 & 15.52 & 22.00 & 13.84 & 39.45 & 81.14 \\
			\end{tabular}
			\vspace{.25 cm}\\
			\caption{Values used in the certified bounds for $K$, $R$, $M$}
			\label{table:3}
		\end{table}
		
		\begin{table}[h]
			\centering
			\renewcommand{\arraystretch}{1}
			\begin{tabular}{c|c}
				$i$ & $\widetilde L_i$ \\
				\hline 
				\\[-.2 cm]
				$0$ & \(\begin{pmatrix}
					-0.2159 & -0.3755 \\
					-0.4694 & 0.4623
				\end{pmatrix}\) \\[.4 cm]
				$1$ & \(\begin{pmatrix}
					1.7718 & -2.1227 \\
					-12.3247 & 16.3690
				\end{pmatrix}\) \\[.4 cm]
				$2$ & \(\begin{pmatrix}
					0.3539 & -4.7730 \\
					0.4185 & -4.6570
				\end{pmatrix}\) \\[.4 cm]
				$3$ & \(\begin{pmatrix}
					-3.1432 & -0.4406 \\
					-0.0916 & -1.1425
				\end{pmatrix}\) \\[.4 cm]
				$4$ & \(\begin{pmatrix}
					-0.4583 & 0.1150 \\
					-0.6163 & 0.8316
				\end{pmatrix}\) \\[.4 cm]
				$5$ & \(\begin{pmatrix}
					1.4772 & -1.9866 \\
					0.3707 & -2.6806
				\end{pmatrix}\) \\[.4 cm]
				$6$ & \(\begin{pmatrix}
					-0.8852 & 0.0258 \\
					-0.1241 & 0.3218
				\end{pmatrix}\) \\[.4 cm]
				$7$ & \(\begin{pmatrix}
					1.0118 & 1.1967 \\
					13.6472 & 13.3154
				\end{pmatrix}\) \\[.4 cm]
				$8$ & \(\begin{pmatrix}
					-0.6239 & -4.3400 \\
					0.7475 & 5.7641
				\end{pmatrix}\) \\[.4 cm]
				$9$ & \(\begin{pmatrix}
					-1.3601 & 1.6743 \\
					-0.3730 & -2.2037
				\end{pmatrix}\) \\
			\end{tabular}
			\vspace{.25 cm}\\
			\caption{The matrices $\widetilde L_i$, chosen to satisfy $\widetilde L_i \approx (Df^c_i)_{0}$}
			\label{table:2}
		\end{table}

	\end{appendices}
	\clearpage
	
	\bibliographystyle{amsplain}
	\bibliography{bibliography_updated}
\end{document}